\documentclass[10pt,reqno,twoside,english]{article}

\usepackage{booktabs}
\usepackage[utf8]{inputenc}
\usepackage[english]{babel}
\usepackage{amsmath,amssymb,amsthm,amsfonts}

\usepackage{adjustbox}

\usepackage[
hmarginratio=1:1,
a4paper,
bottom=3.4cm,
top=3.6cm,
left=2.5cm,
right=2.5cm,
headheight=16pt
]{geometry}

\usepackage{xcolor}
\usepackage{hyperref}
\hypersetup{colorlinks=true, linkcolor=blue}
\usepackage{enumitem}
\usepackage{multirow,rotating,rotfloat}
\usepackage{tikz}
\usepackage{comment}

\usetikzlibrary{positioning,arrows,trees,decorations.text}

\usepackage[square,numbers]{natbib}
\usepackage[toc,page]{appendix}

\newtheorem{theorem}{Theorem}[section]
\newtheorem{lemma}[theorem]{Lemma}
\newtheorem{proposition}[theorem]{Proposition}
\newtheorem{corollary}[theorem]{Corollary}

\theoremstyle{definition}
\newtheorem{definition}[theorem]{Definition}
\newtheorem{example}[theorem]{Example}

\theoremstyle{remark}
\newtheorem{remark}[theorem]{Remark}

\newcommand{\R}{\mathbb{R}}
\newcommand{\calH}{\mathcal{H}}
\newcommand{\calC}{\mathcal{C}}
\newcommand{\calD}{\mathcal{D}}

\newcommand{\calA}{\mathcal{A}}
\newcommand{\Lie}{\mathcal{L}}

\newcommand{\so}{\mathfrak{so}}
\newcommand{\m}{\mathfrak{m}}
\newcommand{\fraku}{\mathfrak{u}}
\newcommand{\f}{\varphi}

\newcommand{\Ker}{\operatorname{Ker}}
\newcommand{\nablag}{\nabla^g}
\newcommand{\tr}{\operatorname{tr}}
\newcommand{\pr}{\operatorname{pr}}
\newcommand{\ddt}{\frac{d}{dt}}

\begin{document}
	
\title{A new approach to the classification of almost contact metric manifolds via intrinsic endomorphisms}
\author{Ilka Agricola, Dario Di Pinto, Giulia Dileo, Marius Kuhrt}
\date{}
\maketitle

%------------------------------------	
\begin{abstract}
	In 1990, D. Chinea and C. Gonzalez gave a classification of almost contact metric manifolds into $2^{12}$ classes, based on the behaviour of the covariant derivative $\nabla^g\Phi$ of the fundamental $2$-form $\Phi$ \cite{CG90}. This large number makes it difficult to deal with this class of manifolds.
 
	We propose a new approach to almost contact metric manifolds by introducing two intrinsic endomorphisms $S$ and $h$, which bear their name from the fact that they are, basically, the entities appearing in the intrinsic torsion. We present a new classification scheme for them by providing a simple flowchart based on algebraic conditions involving  $S$ and $h$, which then naturally leads to a regrouping of the Chinea-Gonzalez classes, and, in each step, to a further refinement, eventually ending in the single classes.
	This method allows a more  natural exposition and derivation of both known and new results, like a new characterization of almost contact metric manifolds admitting a characteristic connection in terms of intrinsic endomorphisms.
	We also describe in detail the remarkable (and still very large) subclass of $\calH$-parallel  almost contact manifolds, defined by the condition $(\nabla^g_X\Phi)(Y,Z)=0$ for all horizontal vector fields,  $X,Y,Z\in\mathcal{H}$. 
\end{abstract}
\medskip

\smallskip


\tableofcontents

%------------------------------------
\section{Introduction}
Almost contact metric manifolds are odd dimensional Riemannian manifolds $(M^{2n+1},g)$ whose structural group is reducible to $U(n)\times 1$. They have been intensively studied in the last decades, after being introduced in the 60's by the Japanese school, as odd dimensional counterparts of almost Hermitian manifolds (see for instance \cite{Sasaki,SH.1} or Blair's monographs \cite{blair76, blair10} and references therein). Some special classes, such as Sasakian, cok\"ahler, or quasi-Sasakian manifolds which include both, have seen a great development---not only for their interplay with K\"ahler geometry, but also as objects of interest in their own, mainly because of their natural occurrence as one-dimensional fibrations, even in the cases where the transverse geometry with respect to the Reeb foliation is not K\"ahler \cite{BoGa08, surveyCDY, Blair67, Kan77, Kan84}. 

One difficulty in dealing with almost contact metric geometry is the huge amount of different classes of manifolds that can be defined, with a great variety of geometric properties. In 1990 a classification into $2^{12}$ classes of almost contact metric manifolds was provided by D. Chinea and C. Gonzalez \cite{CG90}, based on the behaviour of the covariant derivative $\nabla^g\Phi$ of the fundamental $2$-form $\Phi$ associated to the structure $(\varphi,\xi,\eta,g)$, with respect to the Levi-Civita connection. This is an analogue of the Gray-Hervella classification of almost Hermitian manifolds \cite{GH}. In fact, among the $12$ irreducible classes $\calC_1,\ldots,\calC_{12}$ defined by Chinea-Gonzalez, the first four $\mathcal{C}_1,\ldots,\mathcal{C}_4$ correspond to Gray-Hervella's classes $\mathcal{W}_1,\ldots,\mathcal{W}_4$. The reason is that, for the class
\[\mathcal{D}_1:=\mathcal{C}_1\oplus\ldots\oplus\mathcal{C}_4,\]
$\nabla^g \Phi$ is completely determined by the values $(\nabla^g_X\Phi)(Y,Z)$ for any \emph{horizontal} vector fields $X,Y,Z$, namely, vector fields belonging to the distribution $\mathcal{H}=\operatorname{Ker}(\eta)$. On the other hand, the orthogonal class
$$\mathcal{D}_2\oplus\mathcal{D}_3,\quad\text{where } \ 
\mathcal{D}_2= \mathcal{C}_5\oplus\ldots\oplus\mathcal{C}_{11},\quad \mathcal{D}_3= \mathcal{C}_{12},$$
is characterized by
\begin{equation}\label{H-parallel-intro}
(\nabla^g_X\Phi)(Y,Z)=0 \qquad \forall X,Y,Z\in\mathcal{H}.
\end{equation}

Almost contact metric manifolds satisfying \eqref{H-parallel-intro} will be called \emph{$\mathcal{H}$-parallel}, and will be the main class investigated in the present paper. Condition \eqref{H-parallel-intro} is in fact an integrability condition, which coincides with $CR$-integrability in some special cases largely studied in the literature, as it will be discussed in detail in Section \ref{section:main-definition}: this is the case of contact metric, almost cok\"ahler and almost Kenmotsu manifolds.
Unfortunately, the traditional characterizations of the $12$ irreducible classes in terms of $\nabla^g\Phi$ are quite difficult to handle (see Table \ref{Table:CG-classes}),  thus making an alternative approach desirable. 
\\

The first main goal of the present paper is to provide a pair of \emph{intrinsic endomorphisms} $(S,h)$ of the tangent bundle of any almost contact metric manifold, defined by

\[S:=\f\circ\nablag\xi, \qquad h:=\frac12\Lie_\xi\f. \]

The class $\calD_1$ is described by $S=h=0$, and $S$ and $h$ completely determine the covariant derivative $\nabla^g\Phi$ of an $\mathcal{H}$-parallel almost contact metric manifold. Therefore,  they allow to rephrase the characterizations of the corresponding eight irreducible classes in a more handleable way, basically involving purely algebraic properties of the two operators $S$ and $h$ (Table \ref{Table:good_classes}). In fact, we propose a new classification scheme---more precisely, a flowchart---in which the intrinsic endomorphisms $(S, h)$, together with the auxiliary operator
$$P:=-2h+\varphi S\varphi+S, $$
allow to identify step by step mutually orthogonal classes, arriving at the irreducible ones (Figure \ref{Fig:clssif.scheme}, Theorem \ref{theo-flow}). The method even allows to determine the Chinea-Gonzalez class of any almost contact metric manifold, simply by considering the algebraic properties of the $\calD_1$-component of $\nabla^g\Phi$, and those of the intrinsic operators $(S,h)$ which classify the $\calH$-parallel component, namely the component in $\calD_2\oplus\calD_3$.
The same endomorphisms $(S,h)$ make it possible to characterize remarkable geometric properties, as discussed in Section \ref{section:geometric-prop}, namely: $CR$-integrability, normality, anti-normality, integrability of the horizontal distribution $\mathcal{H}$, the characteristic vector field $\xi$ to be Killing.

The reason why we call the endomorphisms $(S,h)$ \emph{intrinsic} is that, together with the structure tensor fields, they completely determine the \emph{intrinsic torsion} of an $\mathcal{H}$-parallel almost contact metric structure.  A similar approach has been adopted in \cite{ACFH}, where the
intrinsic torsions of $SU(3)$-structures in dimension $6$ and $G_2$-structures in dimension $7$ have been expressed in terms of some intrinsic operators. 
 For almost contact metric manifolds, the intrinsic torsion $\Gamma$ is related to the \emph{minimal connection}, a metric connection with torsion adapted to the structure  introduced by J. C. Gonz\'alez-D\'avila and F. Mart\'in Cabrera in \cite{G.MC.}. For each of the eight irreducible classes of $\mathcal{H}$-parallel structures, we determine the algebraic type of the intrinsic torsion, or equivalently, of the torsion of the minimal connection (Theorem \ref{theorem-intrinsic}). 
 We define a new subspace $\calC_{\min}\subset \calC_{10}\oplus\calC_{11}$ 
 that fully describes the
 $\mathcal{H}$-parallel almost contact metric manifolds for which the minimal connection has skew torsion, thus coinciding with the characteristic connection. In fact, using $\calC_{\min}$, we prove that an almost contact metric manifold admits a characteristic connection if and only if it is of class $\calC_1\oplus\calC_3\oplus\calC_4\oplus \calC_{6}\oplus\calC_{7}\oplus \calC_{\min}$ (Theorem \ref{new-char-conn}). This solves the problem, which has been open since the seminal work of Friedrich-Ivanov \cite{FI02}, whether there is a way to characterize such manifolds in terms of Chinea-Gonzalez classes, as it is known for almost Hermitian manifolds or $G_2$-manifolds. 
 It is a remarkable fact that the new class $\calC_{\min}$ is itself irreducible and actually belongs to an infinite family $\calC_{\lambda,\mu}$ of equivalent representations in $\calC_{10}\oplus\calC_{11}$ (Theorem \ref{thm:C_lm}).
 Finally, for $\calH$-parallel almost contact metric manifolds admitting a characteristic connection, we determine a necessary and sufficient condition for the skew torsion to be parallel: this is a condition on the Riemannian curvature, stating that the curvatures $R^g(X,Y)\xi$ are completely determined by the intrinsic operators $(S,h)$ (Theorem \ref{thm:parallel-torsion}).
 
\medskip

%------------------------------------
\subsubsection*{Acknowledgements and declarations.} 
We thank Jenya Ferapontov (Loughborough University, UK) and Arman Taghavi-Chabert (Politechnika {\L}\'odzka, Poland) for an insightful question and fruitful discussions at a conference that led to a refinement of the results of Section \ref{Section:C10+C11}. We also thank the anonymous referee for valuable comments.\\

This material is partially based upon work supported by the National Science Foundation under Grant No. DMS-1928930, while the first author was in residence at the Simons Laufer Mathematical Sciences Institute (formerly MSRI) in Berkeley, California, during the Fall/2024 semester. 
The second author was partially supported by the Centre for Mathematics of the University of Coimbra (funded by the Portuguese Government through FCT project UIDB/00324/2020). The third author is partially supported by the National Center of HPC, Big Data and Quantum Computing, MUR CN00000013, CUP H93C22000450007, Spoke 10. D.D.P. and G.D. are supported by PRIN 2022MWPMAB - \lq\lq Interactions between Geometric Structures and Function Theories\rq\rq, and they are members of INdAM - GNSAGA (Gruppo Nazionale per le Strutture Algebriche, Geometriche e le loro Applicazioni). 

All authors contributed equally to this work and  approved the final manuscript, and we declare to have no competing interests.

%---------------------------------------------------------------------
\section{New tools for the classification of almost contact metric structures}
%----------------------------------------------------------------------
\subsection{Review of the Chinea-Gonzalez classification}

An \emph{almost contact manifold} is a $(2n+1)$-dimensional smooth manifold $M$ endowed with a structure $(\varphi,\xi,\eta)$, where $\varphi$ is a $(1,1)$-tensor field, $\xi$ a vector field, called the \emph{Reeb vector field}, and $\eta$ a $1$-form, such that
\begin{equation*}
	\varphi^2=-I+\eta\otimes \xi,\quad \eta(\xi)=1.
\end{equation*}

It follows that $\varphi\xi =0$, $\eta \circ \varphi =0$, and $\varphi$ has rank $2n$. Moreover, the tangent bundle of $M$ splits as $TM=\calH\oplus\langle\xi\rangle$, where $\calH$ is the $2n$-dimensional distribution defined by $\calH=\mathrm{Im}(\varphi)=\Ker(\eta)$.\\
Any almost contact manifold admits a compatible metric, that is a
Riemannian metric $g$ such that, for every $X,Y\in TM$,
$$g(\varphi X,\varphi Y)=g(X,Y)-\eta(X)\eta(Y).$$
Then, $\eta=g(\cdot,\xi)$ and ${\mathcal H}=\langle \xi\rangle^\perp$. The manifold
$(M,\varphi,\xi,\eta,g)$ is called an \emph{almost contact metric manifold}. The $2$-form
$$\Phi(X,Y):=g(X,\f Y)$$
is called the \textit{fundamental $2$-form} and satisfies $\eta\wedge\Phi^n\neq0$, so that $M$ is orientable. Furthermore, its covariant derivative with respect to the Levi-Civita connection $\nablag$ satisfies the following identities:
\begin{align}
	(\nablag_X\Phi)(Y,Z)&=g(Y,(\nablag_X\f)Z)=-g((\nablag_X\f)Y,Z) \label{eq:C(TM)_1},\\
	(\nablag_X\Phi)(Y,Z)&=-(\nablag_X\Phi)(\f Y,\f Z)+\eta(Y)(\nablag_X\Phi)(\xi,Z)+\eta(Z)(\nablag_X\Phi)(Y,\xi).
	\label{eq:C(TM)_2}
\end{align}

Let us briefly recall the Chinea-Gonzalez classifying criterion of almost contact metric manifolds \cite{CG90}: given a real vector space $V$ of dimension $2n+1$,  endowed with an almost contact metric structure $(\f,\xi,\eta,\langle,\rangle)$, let $\calC(V)$ be the vector space consisting of all tensors $\alpha$ of type $(0,3)$ having the same symmetries as the covariant derivative $\nablag\Phi$ in the case of manifolds, that is
\begin{equation}\label{eq:C(V)}
	\alpha(X,Y,Z)=-\alpha(X,Z,Y)=-\alpha(X,\f Y,\f Z)+\eta(Y)\alpha(X,\xi,Z)+\eta(Z)\alpha(X,Y,\xi)
\end{equation}
for every $X,Y,Z\in V$. The space $\calC(V)$ has a natural inner product induced by that of $V$,
\begin{equation}\label{eq:inner_prod}
	\langle\alpha,\beta\rangle=\sum_{i,j,k=1}^{2n+1}\alpha(e_i,e_j,e_k)\beta(e_i,e_j,e_k),
\end{equation}
where $\alpha,\beta\in\calC(V)$ and $e_1,\dots,e_{2n+1}$ is an arbitrary orthonormal basis of $V$.
Under the natural action of the group $U(n)\times 1$ on $\calC(V)$ defined by
$$(\sigma\alpha)(x,y,z):=\alpha(\sigma^{-1}x,\sigma^{-1}y,\sigma^{-1}z),$$
the space $\calC(V)$ decomposes into twelve orthogonal irreducible factors $\calC_i$ ($i=1,\dots,12$), thus yielding
a total of $2^{12}$ invariant subspaces.\\

For an almost contact metric manifold $(M,\f,\xi,\eta,g)$, this decomposition may be applied to the tangent space in any point $p\in M$. Given an invariant subspace $U$ of $\calC(V)$, we will say that $(M,\f,\xi,\eta,g)$ is of class $U$ if at every point $p\in M$, the covariant derivative $(\nablag\Phi)_p\in\calC(T_pM)$ is of the same type as any $\alpha\in U$. In particular, the null subspace $\calC_0:=\{0\}$ corresponds to the class of cok\"ahler manifolds, defined by $\nablag\Phi=0$. Almost contact metric manifolds corresponding to the irreducible factors $\calC_i$ are characterized by the conditions listed in Table \ref{Table:CG-classes}. In return, these classes of almost contact metric manifolds can be grouped as follows:
\begin{align*}
	\calD_1&=\calC_1\oplus\calC_2\oplus\calC_3\oplus\calC_4,\\
	\calD_2&=\calC_5\oplus\calC_6\oplus\calC_7\oplus\calC_8\oplus\calC_9\oplus\calC_{10}\oplus\calC_{11},\\
	\calD_3&=\calC_{12},
\end{align*}
and the corresponding defining conditions are listed in Table \ref{Table:D-classes}.
In particular, in $\calD_1$ the vertical direction $\xi$ does not contribute to $\nablag\Phi$, so that it can be thought as a $(0,3)$-tensor field on $\calH$.  In fact, the classes $\calC_1$ through $\calC_4$ represent the analog of the Gray-Hervella classes  $\mathcal{W}_1,\dots,\mathcal{W}_4$ which lead to the classification of almost Hermitian manifolds \cite{GH}, and it therefore makes sense to concentrate on the other classes $\calD_2\oplus \calD_3$. This is indeed what we shall do from now on.

\begin{table}[ht]
	\centering
	\setlength{\aboverulesep}{1pt}
	\setlength{\belowrulesep}{1pt}
	\renewcommand{\arraystretch}{1.5}
	\begin{tabular}{clc}
		\toprule
		Class & Defining equation & \\ \midrule
		$\calC_1$ & $(\nabla^{g}_{X}\Phi)(X,Y)=0, \nabla^{g} \eta=0$ & \\ \midrule
		$\calC_2$ & $d \Phi = \nabla^{g} \eta = 0$ & \\ \midrule
		$\calC_3$ & $(\nabla^{g}_{X}\Phi)(Y,Z) - (\nabla^{g}_{\varphi X}\Phi)(\varphi Y, Z)=0,\ \delta \Phi = 0$ & \\ \midrule
		$\calC_4$ & $(\nabla^{g}_{X}\Phi)(Y,Z) = -\frac{1}{2(n-1)} \big[ g(\varphi X,\varphi Y)\delta \Phi(Z) - g(\varphi X, \varphi Z)\delta \Phi(Y)$ & \\
		& $\qquad\qquad\qquad\qquad - \Phi(X,Y)\delta \Phi (\varphi Z) + \Phi(X, Z) \delta \Phi(\varphi Y) \big],
		\ \delta\Phi(\xi)=0$ & \\ \midrule
		$\calC_5$ & $(\nabla^{g}_{X}\Phi)(Y, Z) = \frac{1}{2n}\big[\eta(Y)\Phi(X, Z) - \eta(Z)\Phi(X, Y)\big] \delta \eta$ & \\ \midrule
		$\calC_6$ & $(\nabla^{g}_{X}\Phi)(Y, Z) = \frac{1}{2n}\big[ \eta(Y)g(X, Z) - \eta(Z)g(X, Y) \big] \delta \Phi(\xi)$ & \\ \midrule
		$\calC_7$ & $(\nabla^{g}_{X}\Phi)(Y, Z) = \eta(Z)(\nabla^{g}_{Y}\eta)(\varphi X) + \eta(Y)(\nabla^{g}_{\varphi X}\eta)(Z),\ \delta\Phi=0$ & \\ \midrule
		$\calC_8$ & $(\nabla^{g}_{X}\Phi)(Y,Z) = -\eta(Z)(\nabla^{g}_{Y}\eta)(\varphi X) + \eta(Y) (\nabla^{g}_{\varphi X}\eta)(Z),\ \delta\eta=0$ & \\ \midrule
		$\calC_9$ & $(\nabla^{g}_{X}\Phi)(Y,Z) = \eta(Z)(\nabla^{g}_{Y}\eta)(\varphi X) - \eta(Y)(\nabla^{g}_{\varphi X}\eta)(Z)$ & \\ \midrule
		$\calC_{10}$ & $(\nabla^{g}_{X}\Phi)(Y, Z) = - \eta(Z)(\nabla^{g}_{Y}\eta)(\varphi X)-\eta(Y) (\nabla^{g}_{\varphi X}\eta) (Z)$ & \\ \midrule
		$\calC_{11}$ & $(\nabla^{g}_{X}\Phi)(Y,Z)=- \eta(X)(\nabla^{g}_{\xi}\Phi)(\varphi Y, \varphi Z)$ & \\ \midrule
		$\calC_{12}$ & $(\nabla^{g}_{X}\Phi)(Y, Z) = \eta(X) \eta(Z)(\nabla^{g}_{\xi}\eta)(\varphi Y) -\eta(X)\eta(Y)(\nabla^{g}_{\xi}\eta)(\varphi Z)$ & \\ \bottomrule
	\end{tabular}
	\caption{Chinea-Gonzalez classification of almost contact metric manifolds}
	\label{Table:CG-classes}
\end{table}

\begin{table}[h]
	\centering
	\setlength{\aboverulesep}{1pt}
	\setlength{\belowrulesep}{1pt}
	\renewcommand{\arraystretch}{1.5}
	\begin{tabular}{clc}
		\toprule
		Class & Defining equation & \\ \midrule
		$\calD_1$&
		$(\nablag_\xi\Phi)(Y,Z)=(\nablag_X\Phi)(\xi,Z)=0$ &
		\\ \midrule
		$\calD_2$&
		$(\nablag_X\Phi)(Y,Z)=\eta(X)(\nablag_\xi\Phi)(Y,Z)+
		\eta(Y)(\nablag_X\Phi)(\xi,Z)+\eta(Z)(\nablag_X\Phi)(Y,\xi)$ &
		\\ \midrule
		$\calD_3$&
		$(\nablag_X\Phi)(Y,Z)=\eta(X)\eta(Y)(\nablag_\xi\Phi)(\xi,Z)+
		\eta(X)\eta(Z)(\nablag_\xi\Phi)(Y,\xi)$ &
		\\ \bottomrule
	\end{tabular}
	\caption{Description of the three main classes $\calD_1$, $\calD_2$, $\calD_3$}
	\label{Table:D-classes}
\end{table}

%------------------------------------
\subsection{$\calH$-parallel almost contact metric manifolds} \label{section:main-definition}

We introduce the most important subclass of almost contact metric manifolds of this paper, called \textit{$\calH$-parallel}. Roughly, we impose that $\nablag\Phi$ vanishes in all possible horizontal directions, which is a very natural condition (once we introduced the intrinsic endomorphisms $S$ and $h$ in the following section, we shall prove that it coincides with $\calD_{2}\oplus\calD_{3}$).

\begin{definition}\label{Def:H-parallel}
	An almost contact metric manifold $(M,\f,\xi,\eta,g)$ will be called \textit{$\calH$-parallel} if
	\begin{equation}\label{eq:def.H-par.}
		(\nablag_X\Phi)(Y,Z)=0\quad \forall X,Y,Z\in\calH.
	\end{equation}
	\end{definition}
This is in fact an integrability condition for the structure, which we will relate to the well-known conditions of normality and $CR$-integrability. We will see that in some special cases $\calH$-parallelism coincides with $CR$-integrability.\\

First we will provide an equivalent formulation of $\calH$-parallelism in terms of $d\Phi$ and the tensor field $N_\f$ whose vanishing characterizes normality. Recall that an almost contact manifold $(M,\f,\xi,\eta)$ is said to be \emph{normal} if the almost complex structure $J$, defined on the product manifold $M\times\R$ by
$$J\left(X,f\ddt\right)=\left(\f X-f\xi,\eta(X)\ddt\right),$$
for any $X\in TM$ and $f$ smooth function on $M\times\R$, is  integrable. This is equivalent to the vanishing of the $(1,2)$-tensor field \[N_\f:=[\f,\f]+d\eta\otimes\xi,\] where $[\f,\f]$ is the Nijenhuis torsion of $\f$. Explicitly, for every $X,Y\in TM$,
\begin{equation}\label{eq-normality}
    N_\f(X,Y)=[\f X,\f Y]+\f^2[X,Y]-\f[X,\f Y]-\f[\f X,Y]+ d\eta(X,Y)\xi.
\end{equation}
In particular, for normal structures, $N_\f=0$ implies (\cite[Theorem 6.1]{blair10})
\[\mathcal{L}_\xi\eta=0,\qquad \mathcal{L}_\xi\f=0,\qquad d\eta(\f X,\f Y)=d\eta(X,Y)\quad\forall X,Y\in TM,\]
where $\Lie_\xi$ denotes the Lie derivative with respect to $\xi$.

\begin{proposition}\label{Prop:dPhi=N=0}
	An almost contact metric manifold $(M,\f,\xi,\eta,g)$ is $\calH$-parallel if and only if
	\begin{equation}\label{eq:dPhi=N=0}
		d\Phi(X,Y,Z)=N_\f(X,Y,Z)=0\quad \forall X,Y,Z\in\calH,
	\end{equation}
	where $N_\f(X,Y,Z):=g(N_\f(X,Y),Z)$.
\end{proposition}
\begin{proof}
	Recall that for every almost contact metric manifold the following equation holds for every $X,Y,Z\in TM$ (\cite[Lemma 6.1]{blair10}):
	\begin{align}\label{eq:nablaPhi_acm}
		2g((\nablag_X\f)Y,Z)&=d\Phi(X,\f Y,\f Z)-d\Phi(X,Y,Z)+g(N_\f(Y,Z),\f X)\nonumber\\
		&\quad +\eta(X)d\eta(\f Y,Z)-\eta(X)d\eta(\f Z,Y)\nonumber\\
		&\quad +\eta(Z)d\eta(\f Y,X)-\eta(Y)d\eta(\f Z,X).
	\end{align}
	Thus, for every $X,Y,Z\in\calH$, it reduces to
	$$2(\nablag_X\Phi)(Y,Z)=d\Phi(X,Y,Z)-d\Phi(X,\f Y,\f Z)-N_\f(Y,Z,\f X).$$
	Therefore, if \eqref{eq:dPhi=N=0} holds, then $(\nablag_X\Phi)(Y,Z)=0$ for every $X,Y,Z\in\calH$, and the manifold is $\calH$-parallel. Conversely, for $X,Y,Z\in\calH$, the $\calH$-parallel condition implies that $d\Phi(X,Y,Z)=\mathfrak{S}_{X,Y,Z}(\nablag_X\Phi)(Y,Z)=0$ and hence by the above formula we conclude that $N_\f(Y,Z,\f X)=0$. Being $\f:\calH\to\calH$ an isomorphism, we conclude that $N_\f(X,Y,Z)=0$ for every $X,Y,Z\in\calH$.
\end{proof}
\medskip

$CR$-integrability is a weaker condition than normality. On an almost contact manifold $(M,\f,\xi,\eta)$, the horizontal distribution $\calH$, together with the endomorphism $J_\calH:=\f|_\calH$, which satisfies $J_\calH^2:=-I$, gives rise to an \emph{almost $CR$-structure} $(\calH,J_\calH)$. This is said to be \emph{integrable} if, for all $X,Y\in\calH$,
\begin{itemize}
	\item[(i)] $[J_\calH X,J_\calH Y]-[X,Y]\in\calH$,
	\item[(ii)] $[J_\calH X,J_\calH Y]-[X,Y]-J_\calH ([X,J_\calH Y]+[J_\calH X,Y]) =0$.
\end{itemize}
One also says that $(\calH,J_\calH)$ is a \emph{$CR$-structure}. Next we have a characterization of $CR$-integrability for almost contact manifolds.

\begin{proposition}\label{lemmaCR} 
	Let $(M,\f,\xi,\eta)$ be an almost contact manifold. Then,
	$(\calH,J_\calH)$ is a \emph{$CR$-structure} if and only if \[N_\f(X,Y)=0 \quad \forall X,Y\in\calH.\]
	Furthermore, if $(\calH,J_\calH)$ is a $CR$-structure, then $M$ is normal if and only if $\mathcal{L}_\xi\f=0$.
\end{proposition}
\begin{proof}
     The above condition (i) is equivalent to require $\eta([\f X,\f Y]-[X,Y])=0$, and thus to
	\begin{equation}\label{eq:H_deta-f-invariance}
		d\eta(\f X,\f Y)=d\eta(X,Y)\quad \forall X,Y\in\calH.
	\end{equation}	
	This is also equivalent to $\eta(N_\f(X,Y))=0$ for every $X,Y\in\calH$, by the definition of $N_\f$. On the other hand, if (i) holds, for every $X,Y\in\calH$,
	\begin{align*}
		N_\f(X,Y)&=[\f X,\f Y]+\f^2[X,Y]-\eta([X,Y])\xi-\f[X,\f Y]-\f[\f X,Y]\\
		&=[J_\calH X,J_\calH Y]-[X,Y]-J_\calH ([X,J_\calH Y]+[J_\calH X,Y]),
	\end{align*}
	from which the first statement follows. As regards the second part, assuming the integrability of the almost $CR$-structure, one has that $M$ is normal if and only if $N_\f(\xi,X)=0$ for every $X\in\calH$. By \eqref{eq-normality},
	\begin{equation}\label{eq:N(xi,.)first}
		N_\f(\xi,X)=-\f(\Lie_\xi\f)X+d\eta(\xi,X)\xi
	\end{equation}
	and this vanishes if and only if both the vertical and the horizontal parts vanish, namely $d\eta(\xi,X)=0$ and $\f(\Lie_\xi\f)X=0$. This is equivalent to $\Lie_\xi\f=0$, since one easily check that $d\eta(\xi,X)=\eta((\Lie_\xi\f)\f X)$.
\end{proof}
\medskip

Proposition \ref{lemmaCR} and Proposition \ref{Prop:dPhi=N=0} allow to relate the two conditions of $CR$-integrability and $\calH$-parallelism on three remarkable classes of almost contact metric manifolds largely studied in the literature. Recall that an almost contact metric manifold $(M,\f,\xi,\eta,g)$ is called
\begin{itemize}
	\item \emph{contact metric} if $d\eta=2\Phi$;
	\item \emph{almost cok\"ahler} if $d\eta=0$ and $d\Phi=0$;
	\item \emph{almost Kenmotsu} if $d\eta=0$ and $d\Phi=2\eta\wedge\Phi$.
\end{itemize}

Almost cok\"ahler manifolds are also called almost cosymplectic manifolds.
Here we adopt the definition in \cite{surveyCDY}. In all the above classes, $$d\Phi(X,Y,Z)=0,\qquad  N_\f(X,Y,\xi)=0\quad \forall X,Y,Z\in\calH,$$  
where in particular the second identity is consequence of the definition of $N_\f$ and the fact that equation \eqref{eq:H_deta-f-invariance} is satisfied. As a consequence, we have:

\begin{proposition}
	Let $(M,\f,\xi,\eta,g)$ be an almost contact metric manifold. Assume $M$ to be contact metric, or almost cok\"ahler, or almost Kenmotsu. Then, the structure is $\calH$-parallel if and only if it is $CR$-integrable.
\end{proposition}

In all the three classes, $CR$-integrability has been characterized by means of the covariant derivative $\nablag\varphi$. One should consider the operator
$$h:=\frac12\mathcal{L}_\xi\varphi,$$
which was introduced for contact metric manifolds by Blair in \cite{blair76-h} (se also \cite[Chapter IV, Section 3]{blair76}), and afterwards also adopted for the classes of almost coK\"ahler and almost Kenmotsu manifolds. Then,
\begin{itemize}
	\item a contact metric structure is $CR$-integrable if and only if (\cite{Ta} and \cite[Theorem 6.6]{blair10})
	$$(\nablag_X\varphi)Y=g(X+hX,Y)\xi-\eta(Y)(X+hX)\quad\forall X,Y\in TM;$$
	\item an almost cok\"ahler structure is $CR$-integrable if and only if (\cite{Ol87,surveyCDY})
	$$(\nablag_X\varphi)Y=g(hX,Y)\xi-\eta(Y)hX\quad\forall X,Y\in TM;$$
	\item an almost Kenmotsu structure is $CR$-integrable if and only if (\cite{DP})
	$$(\nablag_X\varphi)Y=g(\f X+hX,Y)\xi-\eta(Y)(\f X+hX)\quad\forall X,Y\in TM.$$
\end{itemize}
In all the above three cases, assuming $CR$-integrability, the structure is normal if and only if $h=0$ (see Proposition \ref{lemmaCR}). This gives rise, respectively, to
\begin{itemize}
	\item \emph{Sasakian} manifolds, i.e.~normal contact metric manifolds, which are characterized by
	$$(\nablag_X\varphi)Y=g(X,Y)\xi-\eta(Y)X\quad\forall X,Y\in TM;$$
	\item \emph{cok\"ahler} manifolds, i.e.~normal almost cok\"ahler manifolds, characterized by $\nablag\f=0$;
	\item \emph{Kenmotsu} manifolds, i.e.~normal almost Kenmotsu manifolds, characterized by
	$$(\nablag_X\varphi)Y=g(\f X,Y)\xi-\eta(Y)\f X\quad\forall X,Y\in TM.$$
\end{itemize}

It is worth remarking that the three classes of $CR$-integrable, or equivalently $\calH$-parallel, contact metric, almost cok\"ahler, almost Kenmotsu manifolds, are quite large, even in the non normal case. Various examples are obtained requiring the Reeb vector field to belong to some nullity distribution for the Riemannian curvature (see for instance \cite{Blairkmi}, \cite{Dacko}, \cite{DP}).
\\

We collect in the following proposition some fundamental properties of the operator $h$ and the covariant derivative $\nabla^g\xi$ for contact metric, almost cok\"ahler, almost Kenmotsu manifolds (see \cite[Lemma 6.2]{blair10}, \cite[Proposition 3.10]{surveyCDY} and \cite{DP}).

\begin{proposition}\label{prop3classes}
	Let $(M,\f,\xi,\eta,g)$ be an almost contact metric manifold. Assume $M$ contact metric, or almost cok\"ahler, or almost Kenmotsu. Then the operator $h$ is symmetric and anticommutes with $\varphi$, namely $h\varphi+\varphi h=0$. The covariant derivative $\nabla^g\xi$ satisfies, respectively, for the three classes
	\begin{equation}\label{3classes}
		\nablag_X\xi=-\varphi X-\varphi hX,\qquad \nablag_X\xi=-\varphi hX,\qquad 	\nablag_X\xi=-\varphi^2 X-\varphi hX.
	\end{equation}
\end{proposition} 

Finally, in the next proposition the vanishing of $h$ (without assuming $CR$-integrability) is characterized.
\begin{proposition}\label{3classes:h=0}
	Let $(M,\f,\xi,\eta,g)$ be an almost contact metric manifold. Then,
    \begin{itemize}
    	\item if $M$ is contact metric, $h=0$ if and only if $\xi$ is Killing (in which case $M$ is called a \emph{$K$-contact} manifold) \cite[Theorem 6.2]{blair10};
    	\item if $M$ is almost cok\"ahler, $h=0$ if an only if $\xi$ is Killing, in which case $M$ is locally isometric to $\mathbb{R}\times N$, with $N$ almost K\"ahler manifolds \cite[Theorem 3.11]{surveyCDY};
     	\item if $M$ is almost Kenmotsu, $h=0$ if an only if $M$ is locally isometric to to a warped product $I\times_f N$, where $N$ is an almost K\"ahler manifold, $I$ is an open interval with coordinate $t$, and $f^2 = ce^{2t}$, for some positive constant $c$ \cite[Theorem 2]{DP07}.
	\end{itemize}
\end{proposition}

It is clear that both the operators $\nablag\xi$ and $h$ have played a central role in the investigation of various classes of almost contact metric manifolds, without receiving a name. They will be fundamental tools in the study of $\calH$-parallel almost contact metric manifolds, as we will see in the next section.

%------------------------------------
\subsection{The intrinsic endomorphisms $(S,h)$}

We introduce the fundamental operators which will be key tools in the classification scheme for almost contact metric manifolds in general and  $\calH$-parallel almost contact metric manifolds in particular.

\begin{definition}
	Let $(M,\f,\xi,\eta,g)$ be an almost contact metric manifold. We set
	\begin{equation}
		S:=\f\circ\nablag\xi, \qquad h:=\frac12\Lie_\xi\f.
	\end{equation}
	We call $(S,h)$ the couple of \textit{intrinsic endomorphisms} of $M$.
\end{definition}

Since $\eta\circ S=0$ and $h\xi=0$, the endomorphisms $S$ and $h$ can be thought as $$S:TM\to\calH,\qquad h:\calH\to TM.$$
They satisfy the following:

\begin{proposition}
	Let $(M,\f,\xi,\eta,g)$ be an almost contact metric manifold. Then, for every $X,Y,Z\in TM$,
	\begin{equation}\label{eq:nablaPhi_S}
		(\nablag_X\Phi)(\xi,Z)=g(SX,Z),
	\end{equation}
	\begin{equation}\label{eq:nabla_xiPhi}
		(\nablag_\xi\Phi)(Y,Z)=-2g(hY,Z)-g(S\varphi Y,\varphi Z)+g(SY,Z).
	\end{equation}
\end{proposition}
\begin{proof}
	The first equation follows by taking $Y=\xi$ in \eqref{eq:C(TM)_1}:
	$$(\nablag_X\Phi)(\xi,Z)=-g((\nablag_X\f)\xi,Z)=g(\f\nablag_X\xi,Z)=g(SX,Z)$$
	Moreover, for every $X,Y\in TM$ one has:
	\begin{align*}
		2hY&=[\xi,\varphi Y]-\varphi[\xi,Y]= (\nablag_\xi\varphi)Y-\nablag_{\varphi Y}\xi+ \varphi\nablag_Y\xi\\
		&= (\nablag_\xi\varphi)Y+\varphi^2\nablag_{\varphi Y}\xi+ \varphi\nablag_Y\xi\\
		&=(\nablag_\xi\varphi)Y+\varphi S\varphi Y+SY,
	\end{align*}
	where in the third equality we used the fact that $\nablag\xi:TM\to \calH$. Therefore, taking the scalar product with $Z\in TM$, one gets \eqref{eq:nabla_xiPhi}.
\end{proof}

In view of the above proposition, it is convenient to introduce a third operator, depending on $S$ and $h$, namely
\begin{equation}\label{P}
	 P:=-2h+\f S\f+S:TM\to TM.
 \end{equation}
 Then, \eqref{eq:nabla_xiPhi} becomes
 \begin{equation}\label{eq:nabla_xiPhi(1)}
	 (\nablag_\xi\Phi)(Y,Z)=g(PY,Z),
\end{equation}
which, in particular, implies that $P$ is skew-symmetric with respect to $g$. From Table \ref{Table:D-classes}, one immediately has:

\begin{proposition}
	An almost contact metric manifold $(M,\f,\xi,\eta,g)$ belongs to the class $\calD_1$ if and only if $S=0$ and $P=0$, which is also equivalent to the vanishing of the intrinsic endomorphisms $S$ and $h$.
\end{proposition}

\begin{remark}\label{S=0h=0} 
	The vanishing of $S$ and $h$ can be interpreted geometrically. Observe that, for any almost contact metric structure, the following hold:
	\begin{equation}\label{eq:S_nabla_eta}
		g(SX,Y)=-(\nablag_X\eta)(\f Y),\qquad g(SX,\f Y)=(\nablag_X\eta)(Y).
	\end{equation}
	Indeed, 
	$$g(SX,Y)=g(\f\nablag_X\xi,Y)=-g(\nablag_X\xi,\f Y)=-X(\eta(\f Y))+ \eta(\nablag_X\f Y)=-(\nablag_X\eta)(\f Y).$$ 
	Then, replacing $Y$ by $\f Y$ and using $(\nablag_X\eta)\xi=0$, one obtains the second identity. Therefore, we have
    \begin{equation}\label{eq:deta}
        d\eta(X,Y)=(\nablag_X\eta)(Y)-(\nablag_Y\eta)(X)=g(SX,\f Y)-g(SY,\f X)=g(X,\f SY)-g(\f SX,Y).
    \end{equation}
	In particular:
	\begin{itemize}
     \item When $S=0$, we have $d\eta=0$. Therefore, the horizontal distribution $\calH$ is integrable, and since $\nablag\xi=0$, the leaves are totally geodesic. The manifold $M$ is then locally isometric to a Riemannian product $\mathbb{R}\times N$, where $N$ is an integral submanifold of $\calH$.  
	\item  The vanishing of $h$, i.e. $\mathcal{L}_\xi\f=0$, expresses the fact that the Reeb vector field $\xi$ is an infinitesimal automorphism of the horizontal distribution $\calH$, and the $(1,1)$-tensor field $\f$ is projectable along the $1$-dimensional foliation spanned by $\xi$. Therefore, the manifold admits local submersions $\pi:U\to U/\xi$, and the base manifold $U/\xi$ is endowed with an almost complex structure, given by the projection of $\f$.
	\item Consequently, any almost contact metric manifold of class $\calD_1$ ($\Leftrightarrow S=h=0$) admits two orthogonal totally geodesic foliations, tangent to $\xi$ and $\calH$, and the almost complex structures induced by $\varphi$ on the horizontal leaves are preserved by the flow of the Reeb vector field $\xi$.
\end{itemize} 
\end{remark}

\begin{remark}\label{remark-d1}
We summarise some  properties of almost contact metric manifolds in the class $\calD_1$ that carry over directly from the corresponding properties for almost Hermitian manifolds in terms of Gray-Hervella classes \cite{Agricola_Srni, FI02, GH}:
\begin{enumerate}
	\item $d\eta=0$ and $\xi$ is a Killing vector field  for any manifold in $\calD_1$ (see the previous remark),
     \item $d\Phi=0$ if and only if the manifold is of class $\calC_2\subset\calD_1$,
     \item $N_\f$ vanishes exactly for $\calC_3\oplus\calC_4\subset\calD_1$,
     \item $N_\f$ defines a $3$-form if and only if the structure is of class $\calC_1\oplus \calC_3\oplus\calC_4\subset\calD_1$.   
\end{enumerate}
\end{remark}

From now on, we will be concerned with the class $\calD_2\oplus\calD_3$, orthogonal to $\calD_1$. This coincides with the class of $\calH$-parallel almost
contact metric manifolds, for which the covariant derivative $\nabla^g\Phi$ is completely determined by the couple $(S,h)$, as stated in the following:

\begin{proposition}\label{Prop:char.D2+D3}
	Let $(M,\f,\xi,\eta,g)$ be an almost contact metric manifold. Then, the following are equivalent:
	\begin{enumerate}
		\item $M$ belongs to the class $\calD_2\oplus\calD_3$;
 		\item $M$ is $\calH$-parallel;
 		\item for every $X,Y,Z\in TM$,
		 \begin{align}\label{eq:nablaPhi}
			(\nabla^g_X\Phi)(Y,Z)&=\eta(X)\left[2g(h\f Y,\f Z) -g(S\f Y,\f Z)+ g(SY,Z)- \eta(Y)g(S\xi,Z)\right] \nonumber\\ 
			&\quad +\eta(Y)g(SX,Z)-\eta(Z)g(SX,Y),
		\end{align}
		which is also equivalent to
		\begin{equation}\label{eq:nablaPhi(1)}
			(\nablag_X\Phi)(Y,Z)=-\eta(X)g(P\f Y,\f Z)+\eta(Y)g(SX,Z)-\eta(Z)g(SX,Y).
		\end{equation}
 	\end{enumerate}
\end{proposition}
\begin{proof}
	From Table \ref{Table:D-classes} one can immediately see that, if $M$ is of class $\calD_2$ or $\calD_3$, then $(\nablag_X\Phi)(Y,Z)=0$ for every $X,Y,Z\in\calH$. Hence the same condition holds for manifolds in $\calD_2\oplus\calD_3$, that is 
	$$\calD_2\oplus\calD_3\subset \calD:=\{\alpha\in\calC(TM)\ |\ \alpha(X,Y,Z)=0\ \forall X,Y,Z\in\calH\},$$ 
	where we denote $\nablag\Phi$ by $\alpha$.  In order to prove the converse, we show that $\calD$ is contained in the orthogonal complement of $\calD_1$ with respect to the inner product \eqref{eq:inner_prod}. This easily follows from the fact that every $\beta\in\calD_1$ vanishes when computed on triplets containing $\xi$. Indeed, for any $\alpha\in\calD$ and $\beta\in\calD_1$, fixed a local orthonormal frame $\{e_i,e_{2n+1}=\xi\}_{i=1,\dots,2n}$, one has
	$$\langle\alpha,\beta\rangle=\sum_{i,j,k=1}^{2n+1}\alpha(e_i,e_j,e_k)\beta(e_i,e_j,e_k)=\sum_{i,j,k=1}^{2n}\alpha(e_i,e_j,e_k)\beta(e_i,e_j,e_k)=0,$$
	being $e_i,e_j,e_k\in\calH$ for every $i,j,k=1,\dots,2n$.

	Now, assume 2. Then, by \eqref{eq:C(TM)_2}, for every $X,Y,Z\in TM$, we have that
	\begin{align*}
		\lefteqn{(\nablag_X\Phi)(Y,Z)}\\
		&=-(\nablag_X\Phi)(\f Y,\f Z)+\eta(Y)(\nablag_X\Phi)(\xi,Z)+\eta(Z)(\nablag_X\Phi)(Y,\xi)\\
		&=-\eta(X)(\nablag_\xi\Phi)(\f Y,\f Z)+\eta(Y)(\nablag_X\Phi)(\xi,Z)+\eta(Z)(\nablag_X\Phi)(Y,\xi)\\
		&=-\eta(X)[-2g(h\f Y,\f Z)-g(S\f^2Y,\f^2Z)+g(S\f Y,\f Z)] +\eta(Y)g(SX,Z)-\eta(Z)g(SX,Y)\\
		&=\eta(X)[2g(h\f Y,\f Z)-g(S\f^2Y,Z)-g(S\f Y,\f Z)] +\eta(Y)g(SX,Z)-\eta(Z)g(SX,Y)\\
		&=\eta(X)[2g(h\f Y,\f Z)+g(SY,Z)-\eta(Y)g(S\xi,Z)-g(S\f Y,\f Z)]\\
		&\quad +\eta(Y)g(SX,Z)-\eta(Z)g(SX,Y),
	\end{align*}
	where in the fourth equality we used the fact that $S:TM\to \calH$. This shows equation \eqref{eq:nablaPhi}, which is equivalent to \eqref{eq:nablaPhi(1)}, by the definition of $P$. The fact that \eqref{eq:nablaPhi} or \eqref{eq:nablaPhi(1)} imply 2 is obvious.
\end{proof}

As an immediate consequence of the given definitions of the operators $S$, $h$, and $P$, we have:

\begin{corollary}\label{Cor:S=h=0}
	An $\calH$-parallel almost contact metric manifold $(M,\f,\xi,\eta,g)$ is cok\"ahler if and only if $S=0$ and $P=0$, or equivalently $S=0$ and $h=0$.
\end{corollary}

The power of the new intrinsic quantities $S$ and $h$ is that they allow to express all relevant geometric quantities for $\calH$-parallel almost contact metric manifolds. Taking into account Proposition \ref{Prop:dPhi=N=0}, the only possibly non-vanishing components of $N_\f$ and $d\Phi$ are expressed in terms of $S,h$, and $P$ in the following result:

\begin{proposition}
	For an $\calH$-parallel almost contact metric manifold, we have the following identities:
	\begin{equation}\label{N1}
   		N_\f(X,Y,\xi)=g((S\f-\f S)X,Y)-g(X,(S\f -\f S)Y),
	\end{equation}
	\begin{equation}\label{N2}
   		N_\f(X,\xi,Y)=2g(\f hX,Y)-\eta(Y)g(S\xi, \f X),
	\end{equation}
	\begin{equation}
    	d\Phi(\xi,X,Y)=g(PX,Y)+g(X,SY)-g(SX,Y).
	\end{equation}
\end{proposition}
\begin{proof}
	Recall that for every almost contact metric manifold, $N_\f$ can be expressed in terms of the covariant derivative by means
	\begin{align}\label{Ncomplete}
		N_\f(X,Y)&=(\nablag_X\f)\f Y-(\nablag_Y\f)\f X+(\nablag_{\f X}\f)Y-(\nablag_{\f Y}\f)X+\eta(X)\nablag_Y\xi- \eta(Y)\nablag_X\xi.
	\end{align}
	Taking the scalar product with $\xi$ we have:
	\begin{align*}
		N_\f(X,Y,\xi)&=(\nablag_X\Phi)(\xi,\f Y)- (\nablag_Y\Phi)(\xi,\f X)+ (\nablag_{\f X}\Phi)(\xi,Y) -(\nablag_{\f Y}\Phi)(\xi,X)\\
		&=g(SX,\f Y)-g(SY,\f X)+g(S\f X,Y)-g(S\f Y,X),				
	\end{align*}
    which gives \eqref{N1}. As regards \eqref{N2}, by \eqref{eq:N(xi,.)first},
	\begin{equation}\label{eq:N(xi,.)}
		N_\f(\xi,X)=-2\f hX+d\eta(\xi,X)\xi.
	\end{equation}
	On the other hand, by \eqref{eq:deta}, $d\eta(\xi,X)=g(S\xi,\f X)$, which implies \eqref{N2}. Finally,
    \begin{align*}
        d\Phi(\xi,X,Y)&=(\nablag_\xi\Phi)(X,Y)+(\nablag_X\Phi)(Y,\xi)+(\nablag_Y\Phi)(\xi,X)\\
        &=g(PX,Y)-g(SX,Y)+g(SY,X). \qedhere
    \end{align*}
\end{proof}

\begin{remark}\label{Rmk:def.P}
	The irreducible classes of $\calH$-parallel almost contact metric manifolds are $\calC_5,\dots,\calC_{12}$. The defining conditions for these classes, listed in Table \ref{Table:CG-classes}, can be rephrased in terms of $S$ and $P$. Using equations \eqref{eq:S_nabla_eta}, one obtains the expressions in Table \ref{Table:CG-classes_bis}.

	\begin{table}[h]
		\centering
		\setlength{\aboverulesep}{1pt}
		\setlength{\belowrulesep}{1pt}
		\renewcommand{\arraystretch}{1.5}
		\begin{tabular}{clc}
			\toprule
			Class & Defining equation & \\ \midrule
			$\calC_5$ & $(\nablag_X\Phi)(Y,Z) = -\frac{1}{2n}\big[\eta(Y)g(\f X,Z)-\eta(Z)g(\f X,Y)\big] \delta\eta$ & \\ \midrule
			$\calC_6$ & $(\nablag_X\Phi)(Y,Z)=\frac{1}{2n}\big[ \eta(Y)g(X,Z)-\eta(Z)g(X,Y) \big]\delta\Phi(\xi)$ & \\ \midrule
			$\calC_7$ & $(\nablag_X\Phi)(Y,Z)=\eta(Y)g(S\f X,\f Z)-\eta(Z)g(SY,X),\ \delta\Phi=0$ & \\ \midrule
			$\calC_8$ & $(\nablag_X\Phi)(Y,Z)=\eta(Y)g(S\f X,\f Z)+\eta(Z)g(SY,X) ,\ \delta\eta=0$ & \\ \midrule
			$\calC_9$ & $(\nablag_X\Phi)(Y,Z)=-\eta(Y)g(S\f X,\f Z)-\eta(Z)g(SY,X) $ & \\ \midrule
			$\calC_{10}$ & $(\nablag_X\Phi)(Y,Z)=-\eta(Y)g(S\f X,\f Z)+\eta(Z)g(SY,X) $ & \\ \midrule
			$\calC_{11}$ & $(\nablag_X\Phi)(Y,Z)=-\eta(X)g(P\f Y,\f Z)$ & \\ \midrule
			$\calC_{12}$ & $(\nablag_X\Phi)(Y,Z)=\eta(X)\eta(Y)g(S\xi,Z)-\eta(X)\eta(Z)g(S\xi,Y) $ & \\ \bottomrule
		\end{tabular}
		\caption{$\calH$-parallel irreducible classes in terms of $S$ and $P$}
		\label{Table:CG-classes_bis}
	\end{table}	
	
	Furthermore, the defining condition for the reducible class $\calD_2=\calC_5\oplus\cdots\oplus\calC_{11}$ in Table \ref{Table:D-classes} rewrites as
	\begin{equation}\label{eq:D_2}
	 (\nablag_X\Phi)(Y,Z)=\eta(X)g(PY,Z)+\eta(Y)g(SX,Z)-\eta(Z)g(SX,Y).
	\end{equation}
\end{remark}

\begin{remark}\label{remark-sum}
	In the following discussion, we will use the fact that the assignment from $\alpha\in\calC(TM)$ to the associated operators $S$, $P$ and $h$ is linear. Indeed, for any tensors $\alpha_1,\alpha_2\in\calC(TM)$, considering the operators $S_i$, $P_i$ and $h_i$, defined by 
	$$g(S_iX,Z)=\alpha_i(X,\xi,Z),\qquad g(P_iY,Z)=\alpha_i(\xi,Y,Z),\qquad P_i=-2h_i+\f S_i\f+S_i,$$
	the operators associated to $\alpha=\alpha_1+\alpha_2$, are obviously given by
	$$S=S_1+S_2,\qquad P=P_1+ P_2,\qquad h=h_1+h_2.$$
	Analogously, if $\alpha=\lambda\alpha_1$, for some $\lambda\in\R$, then $S=\lambda S_1$, $P=\lambda P_1$, $h=\lambda h_1$.
\end{remark}

%------------------------------------
\section{Classification of $\calH$-parallel almost contact metric manifolds}
%------------------------------------
\subsection{The flowchart of $\calH$-parallel structures via intrinsic endomorphisms}

In this section, we will use the operators $S$, $h$, and $P$ to characterize step by step mutually orthogonal subclasses of $\calH$-parallel almost contact metric manifolds. This will eventually lead to the following theorem, which is one of the main results of this paper.
In examples \ref{ex:3ad-Sasaki} and \ref{ex:nS-nc} at the end of this subsection, we will illustrate how one can use the flowchart in practice, whereas the following subsections will deal with some special geometric conditions.

\begin{theorem}\label{theo-flow}
	The intrinsic endomorphisms $(S,h$) classify $\calH$-parallel almost contact metric manifolds according to the flowchart in Figure \ref{Fig:clssif.scheme}.
\end{theorem}

The result is consequence of the following Propositions \ref{scheme1}, \ref{Prop:char.D2}, \ref{scheme3}, \ref{Prop:C5+...+C8//C9+C10}, \ref{scheme5}, \ref{scheme6}, \ref{scheme7}.

\begin{proposition}\label{scheme1}
	Let $(M,\f,\xi,\eta,g)$ be an $\calH$-parallel almost contact metric manifold. Then, $M$ is of class $\calD_3=\calC_{12}$ if and only if $S|_\calH=0$ and $\f h=0$. In this case $h=\frac12g(S\xi,\cdot)\xi$.
\end{proposition}
\begin{proof}
	Comparing \eqref{eq:nablaPhi(1)} with Table \ref{Table:CG-classes_bis}, one has that $M$ belongs to $\calC_{12}$ if and only if
	\begin{equation*}
		-\eta(X)g(P\f Y,\f Z)+\eta(Y)g(SX,Z)-\eta(Z)g(SX,Y) =\eta(X)\eta(Y)g(S\xi,Z)-\eta(X)\eta(Z)g(S\xi,Y),
	\end{equation*}
 	that is
 	$$\eta(X)g(P\f Y,\f Z)+\eta(Y)g(S\f^2X,Z)-\eta(Z)g(S\f^2X,Y)=0$$
	for every $X,Y,Z\in TM$.  Since $\eta\circ S=0$, for $Y=\xi$ one gets $g(S\f^2X,Z)=0$ for all $X,Z\in TM$, that is $S|_\calH=0$. Thus $P|_\calH=-2h$ and for $X=\xi$ the above equation gives $g(h\f Y,\f Z)=0$ for every $Y,Z\in TM$. Being $h\xi=0$, it is equivalent to $\f h=0$. The converse is trivial.\\
	Concerning the last statement, if $\f h=0$, then for every $Y\in TM$, $hY=\eta(hY)\xi$. Thus, taking $Z=\xi$ in \eqref{eq:nabla_xiPhi}, using $S|_\calH=0$, $\eta\circ S=0$ and \eqref{eq:nablaPhi_S}, one has:
	\[2\eta(hY)=2g(hY,\xi)=-(\nablag_\xi\Phi)(Y,\xi)+g(SY,\xi)=(\nablag_\xi\Phi)(\xi,Y)=g(S\xi,Y).\qedhere\]
\end{proof}

\begin{figure}
\vspace{-1cm}
\centering
\begin{tikzpicture}
	[>=stealth, block/.style={rectangle, draw, minimum height=1cm, minimum width=2cm} ]
	\node(1)[block,align=center]at (0,0){$\calD_2\oplus\calD_3$\\ ($\calH$-parallel)};
	\node(2)[block,align=center]at (13,0) {$\calD_3=\calC_{12}$};
	\node(3)[block,align=center]at (0,-3){$\calD_2$};
	\node(4)[block,align=center]at (13,-3){$\calC_{11}$};
	\node(5)[block,align=center]at (0,-12) {$\calC_5\oplus\cdots\oplus\calC_{10}$};
	\node(6)[block,align=center]at (3.5,-7.5) {$\calC_9\oplus\calC_{10}$};
	\node(7)[block,align=center]at (3.5,-16.5) {$\calC_5\oplus\cdots\oplus\calC_8$};
	\node(7.1)[align=center]at (3.5,-17.3){\footnotesize normal $\calH$-parallel};
	\node(8)[block,align=center]at (13,-6){$\calC_{10}$};
	\node(9)[block,align=center]at (13,-9){$\calC_9$};
	\node(9.1)[align=center]at (13,-10){\parbox{2.8cm}{\centering\footnotesize $\calH$-parallel \\almost cok\"ahler}};
	\node(10)[block,align=center]at (7.5,-13.5) {$\calC_5\oplus\calC_8$};
	\node(11)[block,align=center]at (7.5,-19.5) {$\calC_6\oplus\calC_7$};
	\node(11.1)[align=center]at (7.5,-20.3){\footnotesize quasi-Sasakian};
	\node(12)[block,align=center]at (13,-12){$\calC_8$};
	\node(13)[block,align=center]at (13,-15){$\calC_5$};
	\node(13.1)[align=center]at (13,-15.8){\footnotesize $\beta$-Kenmotsu};
	\node(14)[block,align=center]at (13,-18){$\calC_6$};
	\node(14.1)[align=center]at (13,-18.8){\footnotesize $\alpha$-Sasakian};
	\node(15)[block,align=center]at (13,-21){$\calC_7$};
	
	\path[->] (1)edge node[above]{$S|_\calH=0$, $\f h=0$}(2);
	\path[->] (1)edge node[fill=white]{$S\xi=0$}(3);
	\path[->] (3)edge node[above]{$S=0$}(4);
	\path[->] (3)edge node[pos=0.3,fill=white] {\parbox{2.8cm}{\centering $P=0$\\ ($2h=S+\f S\f$)}}(5);
	\draw[->] (5) -| ([xshift=0.5cm] 5.east) |- (6) node[pos=0.25,fill=white]{\parbox{2.8cm}{\centering $S\f+\f S=0$\\$(h=S)$}};
	\draw[->] (5) -| ([xshift=0.5cm] 5.east) |- (7) node[pos=0.25,fill=white]{\parbox{2.8cm}{\centering $S\f-\f S=0$\\ ($h=0$)}};
	\draw[->] (6) -| ([xshift=3cm] 6.east) |- (8) node[pos=0.75,above]{$S\in\Lambda^2\calH^*$};
	\draw[->] (6) -| ([xshift=3cm] 6.east) |- (9) node[pos=0.75,above]{$S\in S^2\calH^*$};
	\draw[->] (7) -| ([xshift=1cm] 7.east) |- (10) node[pos=0.25,fill=white]{$S\in \Lambda^2\calH^*$};
	\draw[->] (7) -| ([xshift=1cm] 7.east) |- (11) node[pos=0.25,fill=white]{$S\in S^2\calH^*$};
	\draw[->] (10) -| ([xshift=0.5cm] 10.east) |- (12) node[pos=0.75,above]{$\tr(\f S)=0$};
	\draw[->] (10) -| ([xshift=0.5cm] 10.east) |- (13) node[pos=0.75,above]{\scriptsize\parbox{2.5cm}{\centering $S=\beta\f$\\$\beta=-\frac1{2n}\tr(\f S)$}};
	\draw[->] (11) -| ([xshift=0.5cm] 11.east) |- (14) node[pos=0.75,above]{\scriptsize\parbox{2.5cm}{\centering $S=\alpha(I-\eta\otimes\xi)$\\$\alpha=\frac1{2n}\tr(S)$}};
	\draw[->] (11) -| ([xshift=0.5cm] 11.east) |- (15) node[pos=0.75,above]{$\tr(S)=0$};
\end{tikzpicture}
\caption{Flowchart for classifying $\calH$-parallel almost contact metric manifolds based on the intrinsic endomorphisms $S$ and $h$ as well as the auxiliary endomorphism $P$.}
\label{Fig:clssif.scheme}
\end{figure}

\begin{proposition}\label{Prop:char.D2}
	Let $(M,\f,\xi,\eta,g)$ be an $\calH$-parallel almost contact metric manifold. Then, $M$ belongs to $\calD_2$ if and only if one of the following equivalent conditions holds:
	\begin{enumerate}
		\item $S\xi=0$;
		\item $\nabla^g_\xi\xi=0$ (i.e. the integral curve of $\xi$ are geodesics);
		\item $\nablag_\xi\eta=0$;
		\item $\xi\lrcorner d\eta=0$;
		\item $\eta\circ h=0$;
		\item $\f h+h\f=0$.
	\end{enumerate}
	Furthermore, in this case $\tr(h)=0$.
\end{proposition}
\begin{proof}
	First we prove the equivalence between conditions 1--6.	Since $\nablag_X\xi\in\calH$ for every $X\in TM$, the definition of $S$ immediately gives that $S\xi=0$ if and only if $\nabla^g_\xi\xi=0$. Moreover,
	$$(\nablag_\xi\eta)X=\xi(\eta(X))-\eta(\nablag_\xi X)=g(\nablag_\xi\xi,X),$$
	and hence 2 and 3 are equivalent. The other equivalences follow from the following identities, which hold for any $X\in TM$:
	\begin{equation}\label{eq:fh+hf,deta,Sxi}
		2(\f h+h\f)X=d\eta(\xi,X)\xi=g(S\xi,\f X)\xi,
	\end{equation}
	\begin{equation}\label{eq:h-fhf}
		2(h-\f h\f)X=g(S\xi,X)\xi,
	\end{equation}
	\begin{equation}\label{eq:eta(hX)}
		2\eta(hX)=g(S\xi,X).
	\end{equation}
	The identities in \eqref{eq:fh+hf,deta,Sxi} are trivially true for $X=\xi$, being $h\xi=0$. For $X\in\calH$ one has:
	\begin{align*}
		2(\f h+h\f)X&=\f[\xi,\f X]-\f^2[\xi,X]+[\xi,\f^2X]-\f[\xi,\f X]\\
		&=[\xi,X]-\eta([\xi,X])\xi-[\xi,X]\\
		&=-\eta([\xi,X])\xi=d\eta(\xi,X)\xi,
	\end{align*}
	and $d\eta(\xi,X)=g(S\xi,\f X)$ by \eqref{eq:deta}. Replacing $X$ by $\f X$ in \eqref{eq:fh+hf,deta,Sxi} one gets \eqref{eq:h-fhf}. Moreover, by definition of $h$, for every $X\in TM$,
	$$2\eta(hX)=\eta([\xi,\f X]-\f[\xi,X])=\eta([\xi,\f X])=-d\eta(\xi,\f X)=g(S\xi,X).$$
	In particular \eqref{eq:fh+hf,deta,Sxi} proves the equivalence between properties 1, 4, 6, and \eqref{eq:eta(hX)} shows that 1 is equivalent to 5.\\
	Finally, comparing equations \eqref{eq:nabla_xiPhi(1)} and \eqref{eq:D_2}, $M$ belongs to $\calD_2$ if and only if $\eta(X)[g(P\f Y,\f Z)+g(PY,Z)]=0$ for every $X,Y,Z\in TM$, that is $P-\f P\f=0$. By the definition of $P$ and $\eta\circ S=0$, this is equivalent to
	$$-2(h-\f h\f)+\eta\otimes S\xi=0.$$
	If this equation holds, evaluating in $\xi$ and using $h\xi=0$, it gives $S\xi=0$. Conversely, if $S\xi=0$, also $h-\f h\f$ vanishes because of \eqref{eq:h-fhf}.\\
	For the last statement, assuming $\dim M=2n+1$ and choosing a local orthonormal frame of type $\{e_i,\f e_i,\xi\}$ ($i=1,\dots,n$), by property 6 we have:
	\begin{align*}
		\mathrm{tr}(h)&=\sum_{i=1}^n[g(he_i,e_i)+g(h\f e_i,\f e_i)]=
		\sum_{i=1}^n[g(he_i,e_i)-g(\f he_i,\f e_i)]\\
		&=\sum_{i=1}^n[g(he_i,e_i)-g(he_i,e_i)]=0.\qedhere
	\end{align*}
\end{proof}

In the following we will focus on the class $\calD_2=\calC_5\oplus\cdots\oplus\calC_{11}$, providing a complete description of the irreducible classes $\calC_5,\dots,\calC_{11}$ in terms of the couple of intrinsic endomorphisms $S$ and $h$. We point out that in $\calD_2$, being $S\xi=P\xi=0$ and $\eta\circ h=\eta\circ P=0$, $S$, $h$ and $P$ can be thought as endomorphisms of $\calH$, i.e. $$S:\calH\to\calH,\qquad h:\calH\to\calH,\qquad P:\calH\to\calH.$$
In particular, $P-\f P\f=0$ is equivalent to $P\f+\f P=0$, that is $P$ anticommutes with $\f$. Furthermore, fixed a local orthonormal frame  $\{e_1,\dots,e_{2n},e_{2n+1}=\xi\}$, for every $\alpha\in \calD_2$ one has that $\alpha(e_i,e_j,e_k)=\alpha(\xi,\xi,e_k)=0$ ($i,j,k=1,\dots,2n$). Then, by using the skew-symmetry in the last two entries, whenever $\alpha$ or $\beta$ is in $\calD_2$, the scalar product \eqref{eq:inner_prod} reduces to
\begin{equation}\label{eq:inn.prod.D2}
	\langle\alpha,\beta\rangle=\sum_{j,k=1}^{2n}\alpha(\xi,e_j,e_k)\beta(\xi,e_j,e_k)+2\sum_{i,k=1}^{2n}\alpha(e_i,\xi,e_k)\beta(e_i,\xi,e_k).
\end{equation}

\begin{proposition}\label{scheme3}
	Let $(M,\f,\xi,\eta,g)$ be an $\calH$-parallel almost contact metric manifold of class $\calD_2$. Then:
	\begin{enumerate}
		\item $M$ belongs to $\calC_{11}$ if and only if $S=0$;
		\item $M$ belongs to $\calC_5\oplus\cdots\oplus\calC_{10}$ if and only if $P=0$ (or equivalently $\nablag_\xi\Phi=0$).
	\end{enumerate}
\end{proposition}
\begin{proof}
	Comparing \eqref{eq:nablaPhi(1)} with Table \ref{Table:CG-classes_bis}, $M$ is of class $\calC_{11}$ if and only if
	$$\eta(Y)g(SX,Z)-\eta(Z)g(SX,Y)=0$$
	for all $X,Y,Z\in TM$.
	Taking $Y=\xi$ and using $\eta\circ S=0$, it follows that $g(SX,Z)=0$ for every $X,Z\in TM$, namely $S=0$. The converse is trivial.\\
	To prove the second statement, it suffices to show that the subspace of $\calD_2$ defined by $P=0$ (i.e. $\nablag_\xi\Phi=0$) is the orthogonal complement of $\calC_{11}$ in $\calD_2$, with respect to the inner product \eqref{eq:inner_prod}. Let 
	$$\calD:=\{\alpha\in\calD_2\ |\ P=0\}=\{\alpha\in\calD_2\ |\ \alpha(\xi,Y,Z)=0\ \forall\, Y,Z\in TM\}.$$
	From Table \ref{Table:CG-classes_bis} (or Table \ref{Table:CG-classes}), one can immediately see that the pure classes $\calC_5,\dots,\calC_{10}$ are contained in $\calD$, and hence so is the sum $\calC_5\oplus\cdots\oplus\calC_{10}$. Conversely, fixed a local orthonormal frame $\{e_1,\dots,e_{2n},e_{2n+1}=\xi\}$, for any $\alpha\in\calD$ and $\beta\in\calC_{11}$, by \eqref{eq:inn.prod.D2} one immediately gets $\langle\alpha,\beta\rangle=0$, since $\alpha(\xi,e_j,e_k)=0$ and $\beta(e_i,\xi,e_k)=0$, being $S=0$ in $\calC_{11}$.
	Therefore $\calD\subset\calC_{11}^\perp=\calC_5\oplus\cdots\oplus\calC_{10}$, where the orthogonal complement is understood within $\calD_2$.
\end{proof}

\begin{remark}
	In $\calC_{11}$, $S=0$ implies $P=-2h$, and hence $h$ is skew-symmetric with respect to $g$. 
\end{remark}
\medskip

\begin{proposition}\label{Prop:C5+...+C8//C9+C10}
	Let $(M,\f,\xi,\eta,g)$ be an $\calH$-parallel almost contact metric manifold of class $\calD_2$, with $P=0$. Then:
	\begin{enumerate}
		\item $M$ belongs to $\calC_5\oplus\cdots\oplus\calC_8$ if and only if $S\f=\f S$, or equivalently  $h=0$;
		\item $M$ belongs to $\calC_9\oplus\calC_{10}$ if and only if $S\f=-\f S$, or equivalently $S=h$.
	\end{enumerate}
	In particular, in $\calC_9\oplus\calC_{10}$, $\tr(S)=0$.
\end{proposition}
	\begin{proof}
	First we point out that $P=0$ means $2h=S+\f S\f$, or equivalently $2h\f=S\f-\f S$, being $S,h:\calH\to\calH$. Therefore $S$ commutes with $\f$ if and only if $h\f=0$, i.e. $h=0$, while $S\f=-\f S$ if and only if $h\f=S\f$, i.e. $h=S$. In the latter case, $\tr(S)=\tr(h)=0$ because of Proposition \ref{Prop:char.D2}. Furthermore, by $P=0$, equation \eqref{eq:nablaPhi(1)} (or \eqref{eq:D_2}) reduces to
	\begin{equation}\label{eq:C5+...C10}
		(\nablag_X\Phi)(Y,Z)=\eta(Y)g(SX,Z)-\eta(Z)g(SX,Y).
	\end{equation}
	On the other hand, according to \cite[p. 26]{CG90}, the defining condition of the class $\calC_9\oplus\calC_{10}$ is
	\begin{align}\label{eq:C9+C10}
		(\nablag_X\Phi)(Y,Z)&=-\eta(Y)(\nablag_{\f X}\Phi)(\xi,\f Z) -\eta(Z)(\nablag_{\f X}\Phi)(\f Y,\xi) \nonumber \\
		&=-\eta(Y)g(S\f X,\f Z)+\eta(Z)g(S\f X,\f Y).
	\end{align}
	Therefore $M$ belongs to $\calC_9\oplus\calC_{10}$ if and only if for every $X,Y,Z\in TM$
	$$\eta(Y)g((S-\f S\f)X,Z)-\eta(Z)g((S-\f S\f)X,Y)=0.$$
	Taking $Y=\xi$, it is equivalent to $S-\f S\f=0$, that is $S\f+\f S=0$.\\
	In order to prove the first item, arguing as in the proof of the previous proposition, let us define $$\calD:=\{\alpha\in\calD_2 \ |\ P=0,\ S\f=\f S\}$$
	and we show that it is the orthogonal complement of $\calC_9\oplus\calC_{10}$ within $\calC_5\oplus\cdots\oplus\calC_{10}$. First, replacing $Y$ by $\xi$ in the defining condition of the classes $\calC_5,\dots,\calC_8$ listed in Table \ref{Table:CG-classes_bis}, one can easily see that $S$ commutes with $\f$, so that $\calC_5\oplus\cdots\oplus\calC_8\subset\calD$.
	Conversely, if $P=0$ and $S\f=\f S$, equation \eqref{eq:C5+...C10} can be written as
	\begin{align*}
		(\nablag_X\Phi)(Y,Z)&=\eta(Y)g(\f SX,\f Z)-\eta(Z)g(\f SX,\f Y)\\
		&=\eta(Y)g(S\f X,\f Z)-\eta(Z)g(S\f X,\f Y),
	\end{align*}
	and it characterizes $\calD$. Now, fixed a local orthonormal frame $\{e_1,\dots,e_{2n},e_{2n+1}=\xi\}$ for any $\alpha\in\calD$ and $\beta\in\calC_9\oplus\calC_{10}$, one has $\alpha(e_i,\xi,e_k)=\alpha(\f e_i,\xi,\f e_k)$ and $\beta(e_i,\xi,e_k)=-\beta(\f e_i,\xi,\f e_k)$ ($i,k=1,\dots,2n$), and hence
	\begin{align*}
		\langle\alpha,\beta\rangle
		=2\sum_{i,k=1}^{2n}\alpha(e_i,\xi,e_k)\beta(e_i,\xi,e_k)
		=-2\sum_{i,k=1}^{2n}\alpha(\f e_i,\xi,\f e_k)
		\beta(\f e_i,\xi,\f e_k)
		=-\langle\alpha,\beta\rangle,
	\end{align*}
	being $\{\f e_i,\xi\}$ another orthonormal frame. Therefore $\langle\alpha,\beta\rangle=0$.
\end{proof}

\begin{proposition}\label{scheme5}
	Let $(M,\f,\xi,\eta,g)$ be an $\calH$-parallel almost contact metric manifold of class $\calC_9\oplus\calC_{10}$. Then:
	\begin{enumerate}
		\item $M$ belongs to $\calC_9$ if and only if $S\in S^2\calH^*$;
		\item $M$ belongs to $\calC_{10}$ if and only if $S\in\Lambda^2\calH^*$,
	\end{enumerate}
	where $S^2\calH^*$ and $\Lambda^2\calH^*$ denote the spaces of symmetric and skew-symmetric endomorphisms of $\calH$, respectively.
\end{proposition}
\begin{proof}
	Comparing equation \eqref{eq:C9+C10} with Table \ref{Table:CG-classes_bis}, $M$ belongs to $\calC_9$ if and only if $\eta(Z)g(S\f X,\f Y)=-\eta(Z)g(SY,X)$ for every $X,Y,Z\in TM$. Since $S\f=-\f S$, it is equivalent to $g(SX,Y)=g(X,SY)$, i.e. $S$ is symmetric with respect to $g$.
	Analogously, one has that $M$ is of class $\calC_{10}$ if and only if $S$ is skew-symmetric with respect to $g$. Indeed, the right hand sides of the defining conditions of $\calC_9$ and $\calC_{10}$ in Table \ref{Table:CG-classes_bis} differ only for the sign of the second term.
\end{proof}

\begin{proposition}\label{scheme6}
	Let $(M,\f,\xi,\eta,g)$ be an $\calH$-parallel almost contact metric manifold of class $\calC_5\oplus\cdots\oplus\calC_8$. Then:
	\begin{enumerate}
		\item $M$ belongs to $\calC_6\oplus\calC_7$ if and only if $S\in S^2\calH^*$;
		\item $M$ belongs to $\calC_5\oplus\calC_8$ if and only if $S\in\Lambda^2\calH^*$.
	\end{enumerate}
\end{proposition}
\begin{proof}
	In the proof of Proposition \ref{Prop:C5+...+C8//C9+C10} we have seen that $\calC_5\oplus\cdots\oplus\calC_8$ is characterized by
	$$(\nablag_X\Phi)(Y,Z)=\eta(Y)g(S\f X,\f Z)-\eta(Z)g(S\f X,\f Y).$$
	On the other hand, it is known that $\calC_6\oplus\calC_7$ is characterized by
	\begin{align}\label{eq:C6+C7}
		(\nablag_X\Phi)(Y,Z)&=\eta(Y)(\nablag_{\f X}\Phi)(\xi,\f Z)+ \eta(Z)(\nablag_Y\Phi)(X,\xi) \nonumber\\
		&=\eta(Y)g(S\f X,\f Z)-\eta(Z)g(SY,X).
	\end{align}
	(see for instance \cite{DiP.D}). Therefore, since $S$ commutes with $\f$, $M$ belongs to $\calC_6\oplus\calC_7$ if and only if $S$ is symmetric with respect to $g$.\\
	Likewise, the subclass $\calD$ of $\calC_5\oplus\cdots\oplus\calC_8$ having $S\in\Lambda^2\calH^*$, is characterized by
	\begin{align}\label{eq:C5+C8}
		(\nablag_X\Phi)(Y,Z)&=\eta(Y)g(S\f X,\f Z)+\eta(Z)g(SY,X)\\
		&=\eta(Y)(\nablag_{\f X}\Phi)(\xi,\f Z)- \eta(Z)(\nablag_Y\Phi)(X,\xi), \nonumber
	\end{align}
	and from Table \ref{Table:CG-classes_bis} it is straightforward to verify that $\calC_5\oplus\calC_8\subset\calD$. Conversely, fixed a local orthonormal frame $\{e_1,\dots,e_{2n},e_{2n+1}=\xi\}$ for any $\alpha\in\calD$ and $\beta\in\calC_6\oplus\calC_7$, one has $\alpha(e_i,\xi,e_k)=-\alpha(e_k,\xi,e_i)$ and $\beta(e_i,\xi,e_k)=\beta(e_i,\xi,e_k)$ ($i,k=1,\dots,2n$), and hence
	\begin{align*}
		\langle\alpha,\beta\rangle
		=2\sum_{i,k=1}^{2n}\alpha(e_i,\xi,e_k)\beta(e_i,\xi,e_k)
		=-2\sum_{i,k=1}^{2n}\alpha(e_k,\xi,e_i)\beta(e_k,\xi,e_i)
		=-\langle\alpha,\beta\rangle.
	\end{align*}
	Therefore $\langle\alpha,\beta\rangle=0$ and thus the equality $\calD=\calC_5\oplus\calC_8$ holds.
\end{proof}

\begin{remark}
	Recall that $\calC_6\oplus\calC_7$ is the class of \emph{quasi-Sasakian} manifolds introduced by D. E. Blair in \cite{Blair67}. They are defined as normal almost contact metric manifolds such that $d\Phi=0$. By the above discussion we have that an almost contact metric manifold $(M,\f,\xi,\eta,g)$ is quasi-Sasakian if and only if it satisfies \eqref{eq:C5+...C10}, or equivalently $$(\nablag_X\f)Y=-\eta(Y)SX+g(SX,Y)\xi,$$
	with $S\in S^2\calH^*$ and such that $S\f=\f S$. This is the    Kanemaki's characterization of quasi-Sasakian manifolds \cite{Kan77,Kan84}.
\end{remark}

It remains to characterize the irreducible classes $\calC_5,\calC_6,\calC_7,\calC_8$. To this aim we prove the following lemma.

\begin{lemma}\label{Lemma:codiff.}
	For every almost contact metric manifold $(M,\f,\xi,\eta,g)$, $\delta\eta=\tr(\f S).$ Moreover, if $M$ is $\calH$-parallel, then  $$\delta\Phi(Z)=\eta(Z)\tr(S)-g(S\xi,Z)$$
	for every $Z\in TM$. In particular $\delta\Phi(\xi)=\tr(S)$.
\end{lemma}
\begin{proof}
	Assuming $\dim M=2n+1$ and fixing a local orthonormal frame $\{e_1,\dots,e_{2n+1}\}$ the codifferential of any differential form $\omega$ is given by $\delta\omega=-\sum_{i=1}^{2n+1}e_i\lrcorner (\nabla^{g}_{e_i}\omega)$. In particular, considering a local orthonormal frame such that $e_{2n+1}=\xi$, one has:
	\begin{align*}
		\delta \eta &= - \sum_{i=1}^{2n+1}(\nablag_{e_i}\eta)(e_i)
		=-\sum_{i=1}^{2n+1} g(e_i,\nablag_{e_i}\xi)
		=\sum_{i=1}^{2n} g(e_i,\f^2\nablag_{e_i}\xi)
		=\sum_{i=1}^{2n} g(e_i, \f Se_i)
		=\tr(\f S),
	\end{align*}
	where in the third equality we use the fact that $\nabla_{e_i}\xi$ is always horizontal. Moreover, if $M$ is $\calH$-parallel, then for every $Z\in TM$,
	\begin{align*}
		\delta\Phi(Z)&=-\sum_{i=1}^{2n+1}(\nablag_{e_i}\Phi)(e_i,Z)
		=-(\nablag_\xi\Phi)(\xi,Z)-\eta(Z)\sum_{i=1}^{2n} (\nablag_{e_i}\Phi)(e_i,\xi)\\
		&=-g(S\xi,Z)+\eta(Z)\sum_{i=1}^{2n}g(Se_i,e_i)
		=-g(S\xi,Z)+\eta(Z)\tr(S),
	\end{align*}
	being $g(S\xi,\xi)=0$.
\end{proof}

\begin{proposition}\label{scheme7}
	Let $(M,\f,\xi,\eta,g)$ be an $\calH$-parallel almost contact metric manifold, with $\dim M=2n+1$. If $M$ is of class $\calC_5\oplus\calC_8$, then:
	\begin{enumerate}
		\item $M$ belongs to $\calC_5$ if and only if $S=\beta\f$, $\beta=-\frac1{2n}\tr(\f S)$;
		\item $M$ belongs to $\calC_8$ if and only if $\tr(\f S)=0$.
	\end{enumerate}
	Furthermore, if $M$ is of class $\calC_6\oplus\calC_7$, then:
	\begin{enumerate}\setcounter{enumi}{2}
		\item $M$ belongs to $\calC_6$ if and only if $S=\alpha(I-\eta\otimes\xi)$, $\alpha=\frac{1}{2n}\tr(S)$;
		\item $M$ belongs to $\calC_7$ if and only if $\tr(S)=0$.
	\end{enumerate}
\end{proposition}
\begin{proof}
	We know that $h=0$ in $\calC_5\oplus\cdots\oplus\calC_8$. Therefore, each subclass is completely characterized by $S$.
	In particular, from Table \ref{Table:CG-classes_bis}, using Lemma \ref{Lemma:codiff.}, the class $\calC_5$ is characterized by $S=\beta\f$, with $\beta=-\frac1{2n}\delta\eta=-\frac1{2n}\tr(\f S)$. The class $\calC_6$ is characterized by $S=\alpha(I-\eta\otimes\xi)$, with $\alpha=\frac1{2n}\delta\Phi(\xi)=\frac1{2n}\tr(S)$.\\
	Finally, comparing equations \eqref{eq:C6+C7} and \eqref{eq:C5+C8} with Table \ref{Table:CG-classes_bis}, one has that, if $M$ is of class $\calC_5\oplus\calC_8$, then it belongs to $\calC_8$ if and only if $\delta\eta=\tr(\f S)=0$. Instead, if $M$ is of class $\calC_6\oplus\calC_7$, then it belongs to $\calC_7$ if and only if $\delta\Phi=0$. By the above lemma, since $S\xi=0$, this is equivalent to $\eta(Z)\tr(S)=0$ for every $Z\in TM$, that is $\tr(S)=0$.
\end{proof}

\begin{remark}
	We recall that $\calC_5$ is the class of \emph{$\beta$-Kenmotsu} manifolds (Kenmotsu for $\beta=1$). They are defined as normal almost contact metric manifolds such that $d\eta=0$ and $d\Phi=2\beta\eta\wedge\Phi$, for a function $\beta$. $\calC_6$ is the class of \emph{$\alpha$-Sasakian} manifolds (Sasakian for $\alpha=1$), defined as normal almost contact metric manifolds such that $d\eta=2\alpha\Phi$, for a function $\alpha$. Finally, $\calC_5\oplus\calC_6$ is the class of \emph{trans-Sasakian} manifolds \cite{Oubina, Marrero}.
\end{remark}

We collect in Table \ref{Table:good_classes} the obtained descriptions in terms of $(S,h)$ of the  irreducible classes of $\calH$-parallel almost contact metric manifolds. The dimensions of the classes have been already computed in \cite{CG90}. We remark that for $n\ge2$, namely $\dim M\ge5$, all the irreducible classes from $\calC_5$ to $\calC_{12}$ are non-trivial. Otherwise, for $n=1$, $\calC_7=\calC_8=\calC_{10}=\calC_{11}=\{0\}$.

\begin{table}[H]
	\centering
	\begin{tabular}{ccc}
		\toprule
		Class & Description & Dimension \\\midrule
		$\calC_5$ & $h=0$, $S=\beta\varphi$, $\beta=-\frac{1}{2n} \mathrm{tr}(\f S)$ & $1$\\\midrule
		$\calC_6$ & $h=0$, $S=\alpha(I-\eta\otimes\xi)$, $\alpha=\frac{1}{2n}\delta \Phi (\xi)$ & $1$\\\midrule
		$\calC_7$ & $h=0$, $S\in S^2\calH^*$, $S\varphi = \varphi S$, $\mathrm{tr}(S)=0$  & $n^2-1$ \\\midrule
		$\calC_8$ & $h=0$, $S\in \Lambda^2\calH^*$, $S\varphi = \varphi S$, $\mathrm{tr}(\varphi S)=0$ & $n^2-1$ \\\midrule
		$\calC_9$ & $h=S$, $S\in S^2\calH^*$, $S\varphi = -\varphi S$ & $n^2+n$ \\\midrule
		$\calC_{10}$ & $h=S$, $S\in\Lambda^2\calH^*$ $S\varphi = -\varphi S$ & $n^2-n$ \\\midrule
		$\calC_{11}$& $S=0$&$n^2-n$\\\midrule
		$\calC_{12}$ & $S|_\calH=0$, $\f h=0$ & $2n$ \\\bottomrule
	\end{tabular}
	\caption{Description of $\calH$-parallel classes in terms of $S$ and $h$.}
	\label{Table:good_classes}
\end{table}

\begin{remark}\label{remark-total-class}
	It is worth remarking that the classification flowchart provides a powerful tool for dealing with any almost contact metric structure, not necessarily $\calH$-parallel. Indeed, let $\alpha=\alpha_1+\alpha_{2}$ be the $(0,3)$- tensor defined by the covariant derivative $\nabla^g\Phi$, where $\alpha_1$ and $\alpha_{2}$ are respectively the $\calD_1$-component and the $\calH$-parallel component, namely the component belonging to $\calD_2\oplus \calD_3$. Then, $\alpha_1$ is classified by the algebraic properties of the tensor $\alpha|_{\calH\times\calH\times\calH}$ according to the characterizations of the $16$ Gray-Hervella classes, while $\alpha_{2}$ is classified by the algebraic properties of the intrinsic tensors $(S,h)$ and the auxiliary tensor $P$.

	For those classes or geometric properties that are described in terms of $N_\f$, note that, due to \eqref{Ncomplete}, for every $X,Y,Z\in\calH$
	\begin{align*}
   	 	N_\f(X,Y,Z)=&{}-\alpha_1(X,\f Y,Z)-\alpha_1(\f X, Y,Z)+\alpha_1(Y,\f X,Z) +\alpha_1(\f Y,X,Z)
	\end{align*}
	which can encode information about the $\calD_1$-component. Recall that, by \eqref{N1} and \eqref{N2}, $N_\f(X,Y,\xi)$ and $N_\f(X,\xi,Y)$ are completely determined by $(S,h)$, which classify the $\calH$-parallel component $\alpha_2$.
\end{remark}
\medskip

We end this section with a few examples, in which one can see our classification scheme in action.

\begin{example}[$3$-$(\alpha,\delta)$-Sasakian manifolds] \label{ex:3ad-Sasaki}

The class of $3$-$(\alpha, \delta)$-Sasaki manifolds has been introduced in \cite{Agr.Dil}. They are defined as almost $3$-contact metric manifolds $(M,\varphi_i,\xi_i,\eta_i,g)$, $i=1,2,3$, satisfying 
\begin{equation}\label{differential_eta}
	d\eta_i=2\alpha\Phi_i+2(\alpha-\delta)\eta_j\wedge\eta_k
\end{equation}
for every even permutation $(i,j,k)$ of $(1,2,3)$, where $\alpha\neq0$ and
$\delta$ are real constants. These manifolds are hypernormal, in the sense that each structure $(\varphi_i,\xi_i,\eta_i,g)$ is normal. For $\alpha=\delta$ one has a $3$-$\alpha$-Sasaki manifold, which is $3$-Sasaki when $\alpha=1$. The three Reeb vector fields $\xi_i$, $i=1,2,3$, of a $3$-$(\alpha, \delta)$-Sasaki manifold are always Killing. We determine the Chinea-Gonzalez class for each structure $(\varphi_i,\xi_i,\eta_i,g)$. First of all notice that in general this is not an $\calH$-parallel structure. Indeed, equation \eqref{differential_eta} implies
$$d\Phi_i=2(\delta-\alpha)(\eta_k\wedge\Phi_j-\eta_j\wedge\Phi_k),$$
and thus $d\Phi_i$ vanishes on triples of vector fields orthogonal to $\xi_i$ if and only if $\alpha=\delta$. Since $N_{\varphi_i}=0$, Proposition \ref{Prop:dPhi=N=0} implies that the structure is $\calH$-parallel only in the $3$-$\alpha$-Sasaki case. In order to determine the $\calD_1$-component of the structure, by \cite[Proposition 2.3.2]{Agr.Dil}, we have that for every vector fields $X,Y,Z$ orthogonal to $\xi_i$, 
$$(\nablag_X\Phi_i)(Y,Z)=-g((\nablag_X\varphi_i)Y,Z)=
2(\alpha-\delta)\,\left[\eta_k(X)g(\varphi_jY,Z)-\eta_j(X)g(\varphi_kY,Z)\right].$$
A direct computation using the fundamental identities
$$\f_k=\f_i\f_j-\eta_j\otimes\xi_i=-\f_j\f_i+\eta_i\otimes\xi_j,\quad 
\xi_k=\f_i\xi_j=-\f_j\xi_i, \quad \eta_k=\eta_i\circ\f_j=-\eta_j\circ\f_i,$$
shows that 
$$(\nablag_{\f_iX}\Phi_i)(\f_iY,Z)=(\nablag_X\Phi_i)(Y,Z),\qquad \delta \Phi (Z)=-\sum_{r=1}^{2n}(\nablag_{e_r}\Phi_i)(e_r,Z)=0$$ 
where $\{e_i\}$, $i=1,\ldots, 2n$, is a horizontal local orthonormal frame containing $\xi_j$, $\xi_k$. Therefore, the $\calD_1$-component belongs to $\calC_3$ (see Table \ref{Table:CG-classes}). The $\calH$-parallel component is classified according the scheme, using the intrinsic operators $S_i$ and $h_i$, together with $P_i$.  We know that
$$\nablag_{\xi_i}\f_i=0,\quad \mathcal{L}_{\xi_i}\f_i=0,\quad \nablag\xi_i=-\alpha\f_i-(\alpha-\delta)(\eta_k\otimes\xi_j-\eta_j\otimes\xi_k),$$
which imply 
$$P_i=0,\quad h_i=0,\quad  S_i=\f_i\circ\nablag\xi_i=-\alpha\f_i^2-(\alpha-\delta)(\eta_j\otimes\xi_j+\eta_k\otimes\xi_k),$$
and in particular $S_i$ is symmetric and commutes with $\f_i$. Following in the classification scheme the subsequent conditions
$$S_i\xi_i=0\rightsquigarrow P_i=0 \rightsquigarrow h_i=0 \rightsquigarrow S_i\in S^2\calH^*$$
we obtain that the $\calH$-parallel component is in $\calC_6\oplus\calC_7$. It is in $\calC_6$ if and only if $\alpha=\delta$. It is in $\calC_7$ if and only if $\tr (S_i)=0$, namely $\delta+2 n\alpha=0$. 

We conclude that each structure $(\varphi_i,\xi_i,\eta_i,g)$ of a $3$-$(\alpha,\delta)$-Sasaki manifold is of class $\calC_3\oplus\calC_6\oplus\calC_7$, which is coherent with the fact that the structure is normal and $\xi_i$ is Killing (see Remark \ref{rem-normality} and Proposition \ref{Prop:Killing} in the next section). Further, the structure has nonzero components in $\calC_3$ and $\calC_7$ when $\alpha\ne\delta$.
\end{example}

\begin{example}[Nearly Sasakian and nearly cosymplectic manifolds] \label{ex:nS-nc}

Nearly Sasakian and nearly cosymplectic manifolds are special classes of almost contact metric manifolds $(M,\varphi,\xi,\eta,g)$, representing odd-dimensional analogues of nearly K\"ahler manifolds. They are defined respectively by
$$(\nablag_X\varphi)Y+(\nablag_Y\varphi)X=2 g(X,Y)\xi-\eta(X)Y-\eta(Y)X,$$
$$(\nablag_X\varphi)Y+(\nablag_Y\varphi)X=0,$$
for every vector fields $X$ and $Y$ on $M$. We determine the Chinea-Gonzalez class of these structures, referring to \cite{Cappel-Dileo, DeN-Dil-Y} for their main properties. In both cases, being $g((\nabla^g_X\varphi)X,Y)=0$ for every $X,Y$ orthogonal to $\xi$, the $\calD_1$-component of the structure belongs to $\calC_1$ (which corresponds to the nearly K\"ahler class in the Gray-Hervella classification). For the $\calH$-parallel component, we analyze now separately the two cases.

For a nearly Sasakian structure the following hold:
$$\nablag_\xi\f=\frac13\Lie_\xi\f=\frac23 h,\qquad \nablag\xi=-\f-\frac23\f h, $$
with $h$ skew-symmetric and anticommuting with $\varphi$.  The vanishing of $h$ is in fact equivalent to require the manifold to be Sasakian, and this is the case in dimension $2n+1>5$, as showed in \cite{DeN-Dil-Y}. Therefore the only non trivial case is $\dim M=5$.
From the above equations, one has
$$P=-\frac23 h, \qquad S=-\f^2+\frac23 h,$$
From $S\xi=0$, we conclude that there is no $\calD_3$-component. In the non-Sasakian case ($h\ne0$), there is a non zero component in $\calC_{11}$. The component in the orthogonal complement of $\calC_{11}$ is completely determined by $S$. Being $-\f^2$ symmetric and commuting with $\varphi$, and $h$ skew-symmetric and anticommuting with $\varphi$, one has nonvanishing components in $\calC_6$ and $\calC_{10}$. Therefore, a nearly Sasakian structure belongs to $\calC_1\oplus\calC_{6}\oplus\calC_{10}\oplus\calC_{11}$.

For a nearly cosymplectic structure the following hold:
$$\nablag_\xi\f=\frac13\Lie_\xi\f=\frac23 h,\qquad \nablag\xi=-\frac23\f h,$$
where again $h$ is skew-symmetric and anticommutes with $\f$. In particular, 
$$P=-\frac23 h, \qquad S=\frac23 h,$$
and arguing as before, we obtain that the structure belongs to $\calC_1\oplus\calC_{10}\oplus\calC_{11}$. In this case the vanishing of $h$ is equivalent to require the structure in $\calC_1$, so that the manifold is locally isometric to the Riemannian product of a real line and a nearly K\"ahler manifold. In dimension $2n+1>5$, nearly cosymplectic manifolds with $h\ne0$ are locally isometric to the Riemannian product $M^5\times N^{2n-4}$, where $M^5$ is a nearly cosymplectic non-cok\"ahler manifold and $N^{2n-4}$ is a nearly K\"ahler manifold \cite{DeN-Dil-Y}.
\end{example}

Further examples will be discussed in the next section.

%----------------------------------
\subsection{Special geometric conditions for  $\calH$-parallel manifolds} \label{section:geometric-prop}

In this section we characterize, in terms of the intrinsic endomorphisms $(S,h)$, some further classes of $\calH$-parallel almost contact metric manifolds, with remarkable geometric properties. In some cases we obtain entire (reducible) classes in the Chinea-Gonzalez classification.

%----------------------------------
\subsubsection{$CR$-integrability}

\begin{proposition}\label{Prop:CR}
	Let $(M,\f,\xi,\eta,g)$ be an $\calH$-parallel almost contact metric manifold. Then, the following conditions are equivalent:
	\begin{enumerate}
		\item $(\calH,J_\calH)$ is a $CR$-structure;
		\item $(S\f-\f S)|_\calH\in S^2\calH^*$;
		\item $M$ belongs to $\calC_5\oplus\cdots\oplus\calC_9\oplus\calC_{11}\oplus\calC_{12}$.
	\end{enumerate}
    All in all the class of $CR$-integrable almost contact metric manifolds is $\calC_3\oplus\calC_4\oplus\calC_5\oplus\cdots\oplus\calC_9\oplus\calC_{11}\oplus\calC_{12}$.
\end{proposition}
\begin{proof}
	Since $M$ is $\calH$-parallel, by Proposition \ref{Prop:dPhi=N=0}, $N_\f(X,Y,Z)=0$ for every $X,Y,Z\in\calH$. Then, by Proposition \ref{lemmaCR}, the almost $CR$-structure $(\calH,J_\calH)$ is integrable if and only if $N_\f(X,Y,\xi)=0$ for every $X,Y\in\calH$. Then, from \eqref{N1}, this is equivalent to 
	$$g((S\f-\f S)X,Y)=g(X,(S\f -\f S)Y)\quad \forall X,Y\in\calH,$$
	i.e. $(S\f-\f S)|_\calH\in S^2\calH^*$.
	From Table \ref{Table:good_classes} (or Figure \ref{Fig:clssif.scheme}), it is easily seen that this condition is fulfilled by all the irreducible classes, except for $\calC_{10}$, which has $S\f=-\f S\in\Lambda^2\calH^*$. Then, denoting by $\calD$ the subclass of $\calH$-parallel almost contact metric manifolds having $S\f-\f S\in S^2\calH^*$, we have that $\calC_5\oplus\cdots\oplus\calC_9\oplus\calC_{11}\oplus\calC_{12}\subset\calD$. To prove the converse, we show that $\calD\subset \calC_{10}^\perp$, where the orthogonal complement is meant within the class $\calD_2\oplus\calD_3$. Let $\{e_i,e_{n+i}=\f e_i,\xi\}$ ($i=1,\dots,n$) be a local orthonormal frame, and let $\alpha\in\calD$, $\beta\in\calC_{10}$. Then,
	$$\alpha(e_i,\xi,\f e_k)+\alpha(\f e_i,\xi,e_k)=\alpha(e_k,\xi,\f e_i)+\alpha(\f e_k,\xi,e_i),$$
	$$\beta(\xi,e_i,e_k)=0,\quad \beta(e_i,\xi,e_k)=-\beta(e_k,\xi,e_i),\quad
	\beta(\f e_i,\xi,e_k)=\beta(e_i,\xi,\f e_k),$$
	and, by equation \eqref{eq:inn.prod.D2}, one has:
	\begin{align*}
		\lefteqn{\langle\alpha,\beta\rangle=2\sum_{i,k=1}^{2n}\alpha(e_i,\xi,e_k)\beta(e_i,\xi,e_k)}\\
		&=2\sum_{i,k=1}^n[\alpha(e_i,\xi,e_k)\beta(e_i,\xi,e_k)+\alpha(\f e_i,\xi,\f e_k)\beta(\f e_i,\xi,\f e_k)\\
		&\quad +\alpha(\f e_i,\xi,e_k)\beta(\f e_i,\xi,e_k)+ \alpha(e_i,\xi,\f e_k)\beta(e_i,\xi,\f e_k)]\\
		&=2\sum_{i,k=1}^n[(\alpha(e_i,\xi,e_k)-\alpha(\f e_i,\xi,\f e_k))\beta(e_i,\xi,e_k)+(\alpha(\f e_i,\xi,e_k)+ \alpha(e_i,\xi,\f e_k))\beta(e_i,\xi,\f e_k)]\\
		&=-2\sum_{i,k=1}^n[(\alpha(e_k,\xi,e_i)-\alpha(\f e_k,\xi,\f e_i))\beta(e_k,\xi,e_i)+(\alpha(\f e_k,\xi,e_i)+ \alpha(e_k,\xi,\f e_i))\beta(e_k,\xi,\f e_i)]\\
		&=-\langle\alpha,\beta\rangle.
	\end{align*}
	Therefore $\langle\alpha,\beta\rangle=0$, and this completes the proof.
    As regards the last statement, removing the $\calH$-parallel assumption, the condition on the $\calD_1$-component for the $CR$-integrability of the structure is that it belongs to $\calC_3\oplus\calC_4$ (see Remarks \ref{remark-d1} and \ref{remark-total-class}).
\end{proof}

\begin{example}[Contact metric, almost cok\"ahler, almost Kenmotsu manifolds]
We already remarked in Section \ref{section:main-definition} that, for the three classes of contact metric, almost cok\"ahler, and almost Kenmotsu manifolds,  $\calH$-parallelism and $CR$-integrability are equivalent conditions. By Proposition \ref{prop3classes}, in the three classes the intrinsic endomorphisms satisfy, respectively,
$$S=-\varphi^2+h,\qquad S=h,\qquad S=\varphi+h,$$
where $h$ is symmetric and anticommutes with $\varphi$. In particular, the above equations imply that $S\xi=0$ and $P=0$, which localizes the structures in the class $\calD_2$. On the other hand, $\varphi$ is skew-symmetric and $\varphi^2$ is symmetric, and both commute with $\varphi$. As a consequence, taking into account the classification scheme, we have:
\begin{itemize}
	\item $\calH$-parallel contact metric manifolds are in the class $(\calC_6\oplus\calC_9)-\calC_9$;
	\item $\calH$-parallel almost cok\"ahler manifolds coincide with the class $\calC_9$;
	\item $\calH$-parallel almost Kenmotsu manifolds are in the class $(\calC_5\oplus\calC_9)-\calC_9$.
\end{itemize}
In all the above cases, removing the $\calH$-parallel assumption, being $d\Phi(X,Y,Z)=0$ for any $X,Y,Z\in\calH$, one has a $\calD_1$-component belonging to $\calC_2$. In particular, it is known that the class of almost cok\"ahler manifolds coincides with the entire class $\calC_2\oplus\calC_9$ \cite{CG90}.
\end{example}

%----------------------------------
\subsubsection{Normal structures}
In Section \ref{section:main-definition} we remarked the fact that the almost $CR$-structure $(\calH,J_\calH)$ associated to a normal almost contact structure is integrable. The converse holds if and only if $h=0$ (Proposition \ref{lemmaCR}). This in fact clarifies the geometric meaning of normality. Indeed, for a normal structure, being $\Lie_\xi\varphi=0$, the manifold admits  local submersions $\pi:U\to U/\xi$, and the base manifold $U/\xi$ is endowed with an almost complex structure (see Remark \ref{S=0h=0}) which turns out to be integrable, because of the $CR$-integrability. Being also $\Lie_\xi d\eta=0$, the $2$-form $d\eta$ is projectable. The fact that $\eta(N_\f(X,Y))=0$ for every $X,Y\in\calH$, or equivalently $d\eta$ is $\f$-invariant on $\calH$ (see equation \eqref{eq:H_deta-f-invariance}), implies that $d\eta$ projects on the complex base manifold onto a closed $2$-form of type $(1,1)$.

\begin{proposition}
	Let $(M,\f,\xi,\eta,g)$ be an $\calH$-parallel almost contact metric manifold. Then the following conditions are equivalent:
	\begin{enumerate}
		\item $M$ is normal;
		\item $h=0$ and $S\f=\f S$;
		\item $M$ is of class $\calC_5\oplus\cdots\oplus\calC_8$.
	\end{enumerate}
\end{proposition}
\begin{proof}
	In view of Proposition \ref{lemmaCR} and Proposition \ref{Prop:CR}, in order to prove the equivalence of 1 and 2, we only need to show that $h=0$ and $(S\f-\f S)|_\calH\in S^2\calH^*$ imply $S\f=\f S$. Since $h=0$, we have $S\xi=0$ by Proposition \ref{Prop:char.D2}. Then, $P=S+\f S\f$ gives $P\f= S\f-\f S$, which is symmetric with respect to $g$. On the other hand, we have already pointed out that in the class $\calD_2$ ($S\xi=0$), $P$ anticommutes with $\f$. Therefore, since $P$ and $\f$ are always skew-symmetric with respect to $g$, so is $P\f$. Hence, it must be $P\f=S\f-\f S=0$.
	Finally, 2 is equivalent to 3 by the discussion in the previous section.	
\end{proof}

\begin{remark}\label{rem-normality}
	It is already known that an almost contact metric manifold is normal if and only if it belongs to $\calC_3\oplus\cdots\oplus\calC_8$ (see \cite{CG90}). Hence the equivalence between 1 and 3 is a direct consequence of the fact that the class of $\calH$-parallel almost contact metric manifolds is equal to $\calC_5\oplus\cdots\oplus\calC_{12}$.
\end{remark}

%----------------------------------
\subsubsection{Anti-normal structures}
The class of anti-normal almost contact structures has been recently defined in \cite{DiP.D}. We say that an almost contact manifold $(M,\f,\xi,\eta)$ is \emph{anti-normal} if $N_\f=2d\eta\otimes\xi$. Differently from the normal case, one can see that such manifolds admit local submersions $\pi:U\to U/\xi$, with fibers tangent to $\xi$, where the base manifold is a complex manifold endowed with a closed $2$-form of type $(2,0)$.

\begin{proposition}
	Let $(M,\f,\xi,\eta,g)$ be an $\calH$-parallel almost contact metric manifold. Then the following conditions are equivalent
	\begin{enumerate}
		\item $M$ is anti-normal;
		\item $h=0$ and $S\in\Lambda^2\calH^*$.
	\end{enumerate}		
\end{proposition}
\begin{proof}
	The anti-normal condition is equivalent to
	$$N_\f(\xi,X)=0,\qquad N_\f(X,Y)=2d\eta(X,Y)\xi\quad \forall X,Y\in\calH$$
	(see \cite{DiP.D}). As in the proof of Proposition \ref{lemmaCR}, the first equation is equivalent to $h=0$, in which case $P=S+\f S\f$. Owing to Proposition \ref{Prop:dPhi=N=0}, for $\calH$-parallel almost contact metric manifolds, the second condition is equivalent to
	$$N_\f(X,Y,\xi)=2d\eta(X,Y)\quad \forall X,Y\in\calH,$$
	or also to $d\eta(\f X,Y)-d\eta(X,\f Y)=0$, which directly follows form the definition of $N_\f$. Using \eqref{eq:deta},
	\begin{align*}
		 d\eta(\f X,Y)-d\eta(X,\f Y)
		&=g(S\f X,\f Y)-g(\f^2X, SY)-g(SX,\f^2Y)+g(\f X, S\f Y)\\
		&=g((S-\f S\f)X,Y)+g(X,(S-\f S\f)Y),
	\end{align*}
	which vanishes if and only if $P':=S-\f S\f\in\Lambda^2\calH^*$. Since $P$ is always skew-symmetric and $P+P'=2S$, we have that $N_\f(X,Y,\xi)=0$ if and only if $S\in\Lambda^2\calH^*$.
\end{proof}

Notice that, if an almost contact manifold is both normal and anti-normal, necessarily $d\eta=0$. We analyze this condition in Section \ref{deta=0}.

%----------------------------------
\subsubsection{Structures with Killing vector field $\xi$}
In \cite{DiP.D} the authors considered the class of \emph{generalized quasi-Sasakian} manifolds, extending to any dimension the notion introduced in dimension $5$ by C. Puhle \cite{Pu13}. A generalized quasi-Sasakian manifold is defined as an almost contact metric manifold $(M,\f,\xi,\eta,g)$ with Killing Reeb vector field $\xi$, and such that $d\Phi(X,Y,Z)=N_\f(X,Y,Z)=0$ for every $X,Y,Z\in\calH$. Owing to Proposition \ref{Prop:dPhi=N=0}, generalized quasi-Sasakian manifolds coincide with the class of $\calH$-parallel almost contact metric manifolds, with $\xi$ Killing. It is already known that they coincide with the class $\calC_6\oplus\calC_7\oplus\calC_{10}\oplus\calC_{11}$. The following proposition provides a characterization in terms of the intrinsic endomorphism $S$.

\begin{proposition}\label{Prop:Killing}
	Let $(M,\f,\xi,\eta,g)$ be an $\calH$-parallel almost contact metric manifold. Then the following conditions are equivalent:
	\begin{enumerate}
		\item $\xi$ is Killing;
		\item $S\xi=0$ and $\f S\in\Lambda^2\calH^*$ (or equivalently $\f S\in\Lambda^2TM^*$);
		\item $S\xi=0$ and $S\f\in\Lambda^2\calH^*$;
		\item $M$ is of class $\calC_6\oplus\calC_7\oplus\calC_{10}\oplus\calC_{11}$.
	\end{enumerate}
	If one of the above conditions holds, then $h$ and $\f h=-h\f$ are skew-symmetric with respect to $g$.

    All in all, the class of almost contact metric manifolds with $\xi$ Killing is $\calD_1\oplus\calC_6\oplus\calC_7\oplus\calC_{10}\oplus\calC_{11}$.
\end{proposition}
\begin{proof}
    The Reeb vector field $\xi$ is Killing if and only if $\nablag\xi$ is skew-symmetric with respect to $g$, that is, for every $X,Y\in TM$,
	$$0=g(\nablag_X\xi,Y)+g(X,\nablag_Y\xi)=g(\f\nablag_X\xi,\f Y)+g(\f X,\f\nablag_Y\xi)=g(SX,\f Y)+g(\f X,SY),$$
	where in the second equality we used that $\nablag_X\xi$ is always horizontal. Therefore, $\xi$ is Killing if and only if $\f S$ is skew-symmetric with respect to $g$. In particular, for $X=\xi$ in the above equation, one immediately gets $S\xi=0$, so that $S$ can be thought as an endomorphism of $\calH$. Moreover,  replacing $X,Y$ by $\f X,\f Y$ in $g(\f SX,Y)+g(X,\f SY)=0$, owing to $\eta\circ S=0$, one gets 3:
	$$0=g(\f S\f X,\f Y)+g(\f X,\f S\f Y)=g(S\f X,Y)+g(X, S\f Y).$$
	The converse is analogous and it makes use of $S\xi=0$.
	Now, it is known that the class $\calC_6\oplus\calC_7\oplus\calC_{10}\oplus\calC_{11}$ is characterized by
	\begin{align*}
		(\nablag_X\Phi)(Y,Z)&=-\eta(X)(\nablag_\xi\Phi)(\f Y,\f Z)-\eta(Y)(\nablag_{\f Z}\Phi)(\f X,\xi)-\eta(Z)(\nablag_X\Phi)(\xi,Y)\\
		&=-\eta(X)g(P\f Y,\f Z)+\eta(Y)g(S\f Z,\f X)-\eta(Z)g(SX,Y),
	\end{align*}
	(see for instance \cite{DiP.D}).
	Comparing with equation \eqref{eq:nablaPhi(1)}, we get that an $\calH$-parallel almost contact metric manifold belongs to $\calC_6\oplus\calC_7\oplus\calC_{10}\oplus\calC_{11}$ if and only if $$\eta(Y)g(S\f Z,\f X)=\eta(Y)g(SX,Z)$$
	for every $X,Y,Z\in TM$. For $X=Y=\xi$ one gets $S\xi=0$. Moreover, replacing $Y$ with $\xi$ and $X$ with $\f X$, it is equivalent to 
	$$-g(S\f Z,X)=g(S\f X,Z),$$
	i.e. $S\f\in\Lambda^2\calH^*$.

	Concerning the last statement, direct computations show that for any almost contact metric structure
	\begin{equation}\label{eq:dPhi-Lg-h}
		d\Phi(\xi,Y,Z)=(\Lie_\xi g)(Y,\f Z)+2g(Y,hZ)\quad \forall Y,Z\in TM.
	\end{equation}
	Therefore, if $\xi$ is Killing, $h$ is skew-symmetric. On the other hand, $h$ and $\f$  anticommute each other (by Proposition \ref{Prop:char.D2}), and both of them are skew-symmetric with respect to $g$. Then $\f h=-\f h$ is skew-symmetric.
	The last statement is an immediate consequence of the fact that $\xi$ is always Killing for structures of class $\calD_1$ (Remark \ref{remark-d1}).
\end{proof}
\medskip

From equation \eqref{eq:dPhi-Lg-h}, it follows that, for a generalized-quasi-Sasakian manifold, the fundamental 2-form $\Phi$ is closed if and only if $h=0$. In \cite{DiP.D}, it is proved that the class of generalized-quasi-Sasakian manifolds with $h=0$ is exactly the class of \emph{transversely K\"ahler} manifolds.
These are defined as almost contact metric manifolds $(M,\f,\xi,\eta,g)$ such that the structure tensor fields $(\f,g)$ are projectable along the $1$-dimensional foliation generated by $\xi$ and, locally, induce a transverse K\"ahler structure on the space of leaves $M/\xi$. Within the class of transversely K\"ahler manifolds, the normality and anti-normality conditions (or equivalently the $\f$-invariance and the $\f$-anti-invariance of $d\eta$) distinguish the classes of quasi-Sasakian and anti-quasi-Sasakian manifolds, respectively \cite{DiP.D,DiP}. Precisely, \emph{anti-quasi-Sasakian} manifolds coincide with manifolds of class  $\calC_{10}\oplus\calC_{11}$ with $h=0$.

We provide now a characterization of the class $\calC_{10}\oplus\calC_{11}$. It plays a crucial role for the investigation of adapted connections in the next section.

\begin{proposition}\label{Prop:C10+C11}
	An $\calH$-parallel almost contact metric manifold $(M,\f,\xi,\eta,g)$ belongs to $\calC_{10}\oplus\calC_{11}$ if and only if 
	$S\in\Lambda^2\calH^*$ and one of the two following equivalent conditions holds,
	\begin{enumerate}
		\item $S\f+\f S=0$,
		\item $P=-2h+2S$.
	\end{enumerate}
	Furthermore, the following subset of $\calD_2\oplus\calD_3$ defines a subspace, and it lies in $\calC_{10}\oplus\calC_{11}$:
	$$\calC_{\min} := \{ \alpha\in \calD_{2}\oplus\calD_{3} \, | \, h=2S\} \subset \calC_{10}\oplus\calC_{11},$$
	and the defining condition $h=2S$ can equivalently be reformulated as $P=-2S$.
\end{proposition}
\begin{proof} 
	The class $\calC_{10}\oplus\calC_{11}$ is characterized by
	\begin{align*}
		(\nablag_X\Phi)(Y,Z)&=-\eta(X)(\nablag_\xi\Phi)(\f Y,\f Z)-\eta(Y)(\nablag_{\f X}\Phi)(\xi,\f Z)+\eta(Z)(\nablag_Y\Phi)(\xi,X)\\
		&=-\eta(X)g(P\f Y,\f Z)-\eta(Y)g(S\f X,\f Z)+\eta(Z)g(SY,X)
	\end{align*}
	(see \cite{DiP.D}). Comparing with \eqref{eq:nablaPhi(1)}, one has that an $\calH$-parallel almost contact metric manifold $M$ belongs to $\calC_{10}\oplus\calC_{11}$ if and only if
	$$-\eta(Y)g(S\f X,\f Z)+\eta(Z)g(SY,X) =\eta(Y)g(SX,Z)-\eta(Z)g(SX,Y)$$
	for every $X,Y,Z\in TM$.
	Replacing $Y$ by $\xi$ and $Z$ by $\f Z$, one has:
	$$-g(S\f X,\f^2Z)=g(SX,\f Z),$$
	or equivalently,
	$$g(S\f X,Z)=-g(\f SX,Z),$$
	namely $S$ anticommutes with $\f$. Instead, for $Z=\xi$ in the same equation, one obtains $g(SY,X)=-g(SX,Y)$, i.e. $S\in\Lambda^2\calH^*$. The converse is trivial.

    For the last statement, let us first  show that any structure satisfying $h=2S$ has to be of type $\calC_{10}\oplus\calC_{11}$. If $h=2S$,  $S\xi=0$ (because $h\xi=0$), and therefore Proposition \ref{Prop:char.D2} implies $h\f+\f h=0$, which means of course $S\f+\f S=0$. Furthermore, the expression for $P$ may be simplified,
    $$P=-2h+S+\f S \f=-2h+2S=-2S.$$
    As $P$ is always skew-symmetric, this proves that $S$ is skew-symmetric, too.  By the first part, this proves that such a structure is indeed of type $\calC_{10}\oplus\calC_{11}$. The fact that $\calC_{\min}$ is indeed a subspace is immediate consequence of the characterizing condition $P=-2S$ (see Remark \ref{remark-sum}).
\end{proof}

\begin{remark}
	According to Proposition \ref{Prop:C10+C11}, almost contact metric manifolds in $\calC_{10}\oplus\calC_{11}$ having $h=0$ are anti-normal.
	These are the already mentioned anti-quasi-Sasakian structures. Notice that, if they are of pure type $\calC_{10}$ or $\calC_{11}$, then they are cok\"ahler, because in this case they would have $S=h=0$; furthermore, the space of anti-quasi-Sasakian structures intersects the space $\calC_{\min} $ trivially by definition.
\end{remark}
\smallskip

\begin{remark}
 A $K$-contact manifold is a contact metric manifold $(M,\f,\xi,\eta,g)$ with Killing vector field $\xi$. For a contact metric manifold, requiring $\xi$ Killing is equivalent to $h=0$. This has been recalled in Proposition \ref{3classes:h=0}, and actually follows from equation \eqref{eq:dPhi-Lg-h}, being $d\Phi=0$. On the other hand, an $\calH$-parallel contact metric manifold with $h=0$ is necessarily Sasakian (see Section \ref{section:main-definition}). Therefore, an $\calH$-parallel $K$-contact manifold is Sasakian.
\end{remark}

\subsubsection{Infinitely many equivalent representations in $\calC_{10}\oplus\calC_{11}$}\label{Section:C10+C11}

The class $\calC_{\min}$ defined in the above section will play a central role in the description of almost contact metric structures admitting a characteristic connection (Theorem \ref{new-char-conn}). 
In fact, $\calC_{10}$ and $\calC_{11}$
are equivalent as $U(n)\times1$--modules, and hence 
their realisation in the larger class $\calD_2$ is not unique, as it is always the case when the decomposition of a representation into irreducible submodules is not multiplicity-free. Rather, we will now describe a $2$-parameter family of equivariant embeddings of the abstract irreducible representation $\calC_{10}$ into $\calC_{10}\oplus\calC_{11}$ that includes $\calC_{\min}$ as a particular case. \\

\begin{theorem}\label{thm:C_lm}
    For every pair $(0,0)\neq (\lambda,\mu)\in\R^2$ there exists a $U(n)\times 1$-equivariant embedding
    $f_{\lambda,\mu}:\calC_{10}\to\calC_{10}\oplus\calC_{11}$.
    Its image coincides with 
    \begin{equation}\label{eq:C_lm}
        \calC_{\lambda,\mu}=\{\alpha\in{\calD_2\oplus\calD_3}\ |\ \lambda P=\mu S\},
    \end{equation}
    which is then an invariant irreducible subspace of $\calC_{10}\oplus\calC_{11}$ isomorphic to $\calC_{10}$. 
    The orthogonal complement of $\calC_{\lambda,\mu}$ in $\calC_{10}\oplus\calC_{11}$ is 
    $$\calC_{\lambda,\mu}^\perp =\calC_{\lambda',\mu'}\quad  \text{with}\quad 2\lambda\lambda'+\mu\mu'=0$$ and hence it is again isomorphic to $\calC_{10}$. Consequently
$\calC_{10}\oplus\calC_{11}=\calC_{\lambda,\mu}\oplus \calC_{\lambda',\mu'}$ for every coefficients satisfying $2\lambda\lambda'+\mu\mu'=0$.
\end{theorem}

\begin{proof}
For every $(\lambda,\mu)\neq (0,0)$ we define the linear map
$$f_{\lambda,\mu}:\calC_{10}\longrightarrow
\calC_{10}\oplus\calC_{11}\qquad \alpha\longmapsto
\alpha_{\lambda,\mu},$$
where 
$$\alpha_{\lambda,\mu}(X,Y,Z) :=\lambda(\eta(Z)\alpha(X,Y,\xi)+\eta(Y)\alpha(X,\xi,Z))+\mu\eta(X)\alpha(Y,\xi,Z)$$
for any $X,Y,Z\in TM$. In particular
\begin{equation}\label{eq:alpha_lm}
    \alpha_{\lambda,\mu}(\xi,Y,Z)=\mu\alpha(Y,\xi,Z),\qquad 
    \alpha_{\lambda,\mu}(X,\xi,Z)=\lambda\alpha(X,\xi,Z).
\end{equation}
We show that $f_{\lambda,\mu}$ preserves the inner product, up to a non-zero coefficient, so that $f_{\lambda,\mu}$ is injective. Indeed, by the expression of the inner product on $\calD_2$ (see equation \eqref{eq:inn.prod.D2}), for any $\alpha,\beta\in\calC_{10}$ one has:
\begin{align*}
    \langle \alpha_{\lambda,\mu},\beta_{\lambda,\mu}\rangle&= \sum_{i,j=1}^{2n}\alpha_{\lambda,\mu}(\xi,e_i,e_j)\beta_{\lambda,\mu}(\xi,e_i,e_j) + 2\sum_{i,j=1}^{2n}\alpha_{\lambda,\mu}(e_i,\xi,e_j)\beta_{\lambda,\mu}(e_i,\xi,e_j)\\
    &=\mu^2\sum_{i,j=1}^{2n}\alpha(e_i,\xi,e_j)\beta(e_i,\xi,e_j) + 2\lambda^2\sum_{i,j=1}^{2n}\alpha(e_i,\xi,e_j)\beta(e_i,\xi,e_j)\\
    &=(\mu^2+2\lambda^2)\sum_{i,j=1}^{2n}\alpha(e_i,\xi,e_j)\beta(e_i,\xi,e_j)
    \\ &=\frac{\mu^2+2\lambda^2}{2}\langle\alpha,\beta\rangle,
\end{align*}
where the last equality follows from the fact that $\alpha(\xi,e_i,e_j)=0$ for $\alpha\in\calC_{10}$ (i.e. $P=0$). The map $f_{\lambda,\mu}$ is clearly equivariant with respect to the action of $U(n)\times1$. In particular, $f_{\lambda,\mu}(\calC_{10})$ is an invariant irreducible subspace of $\calC_{10}\oplus\calC_{11}$ isomorphic to $\calC_{10}$.
 
We prove that  $f_{\lambda,\mu}(\calC_{10})$ coincides with $\calC_{\lambda,\mu}$ defined in \eqref{eq:C_lm}. 
Indeed, by equations \eqref{eq:alpha_lm} every $\alpha_{\lambda,\mu}\in f_{\lambda,\mu}(\calC_{10})$ satisfies 
$$\lambda\alpha_{\lambda,\mu}(\xi,Y,Z)=\mu\alpha_{\lambda,\mu}(Y,\xi,Z),$$ namely $\lambda P=\mu S$. On the other hand, we show that $\calC_{\lambda,\mu}$ is contained in $\calC_{10}\oplus\calC_{11}$ and has the same dimension as $\calC_{10}$, which is enough to conclude that $f_{\lambda,\mu}(\calC_{10})=\calC_{\lambda,\mu}$. If $\mu=0$ then $P=0$, and thus $\calC_{\lambda,0}=\calC_{10}$; if $\lambda=0$, then $S=0$ and $\calC_{0,\mu}=\calC_{11}$. Otherwise, $S$ is skew-symmetric since so is $P$. Therefore, for every $X\in TM$, $g(S\xi,X)=-g(\xi,SX)=0$, so that $S\xi=0$ and $h\f+\f h=0$ by Proposition \ref{Prop:char.D2}. Furthermore by equation \eqref{P}, $\mu S=\lambda P=\lambda(-2h+S+\f S\f)$, which gives:
$$\mu S\f=\lambda(-2h\f+S\f+\f S\f ^2)=\lambda(-2h\f +S\f-\f S),$$
$$\mu \f S=\lambda(-2\f h+\f S+\f^2S\f)=\lambda(-2\f h+\f S-S\f).$$
By adding these two equations, one gets $\mu(S\f+\f S)=0$, namely $S$ anticommutes with $\f$. Therefore, $\calC_{\lambda,\mu}\subset \calC_{10}\oplus\calC_{11}$, according to Proposition \ref{Prop:C10+C11}. As regards the dimension of $\calC_{\lambda,\mu}$, let $\{e_1,\dots, e_{2n+1}\}$ be an orthonormal frame, with $e_{n+i}=\f e_i$ for $i=1,\dots,n$ and $e_{2n+1}=\xi$. For $\alpha\in\calC_{\lambda,\mu}$ the condition $\lambda P=\mu S$ gives 
$$\lambda\alpha(\xi,e_i,e_j)=\mu\alpha(e_i,\xi,e_j).$$
Moreover, being $S\xi=0$, $S\in\Lambda^2\calH^*$ and $S\f+\f S=0$, it turns out that $\alpha$ is completely determined by $\alpha(e_i,e_j,\xi)$ and $\alpha(e_i,\f e_j,\xi)$ for $1\le i<j\le n$. Thus, 
$$\dim\calC_{\lambda,\mu}=2\frac{(n-1)n}{2}=n^2-n=\dim\calC_{10}.$$
Finally, straightforward computations show that the orthogonal complement of a class $\calC_{\lambda,\mu}$ inside $\calC_{10}\oplus\calC_{11}$ is $\calC_{\lambda,\mu}^\perp=\calC_{\lambda',\mu'}$ with $$2\lambda\lambda'+\mu\mu'=0,$$
so that $\calC_{\lambda,\mu}^\perp$ is again irreducible and isomorphic to $\calC_{10}$.
\end{proof}

\begin{remark}
    Within the infinite family of equivalent irreducible representations in $\calC_{10}\oplus\calC_{11}$,  $$\calC_{1,-2}=\calC_{\min},\qquad \calC_{1,2}=\{\text{anti-quasi-Sasakian structures}\}.$$
    In Figure \ref{fig:C_lm} three couples of orthogonal classes of type $\calC_{\lambda,\mu}$ are displayed. Notice that, as computed in Example \ref{ex:nS-nc}, for both nearly-Sasakian and nearly cosymplectic structures the component in $\calC_{10}\oplus\calC_{11}$ lies in $\calC_{1,-1}$, which is the orthogonal class of anti-quasi-Sasakian structures. 
\end{remark}

\begin{figure}
\centering
\begin{tikzpicture}[domain=-3:3,scale=0.7]
    \draw[->] (-4,0) -- (4,0) node[below right] {\footnotesize $S$};
    \node at(-4.5,0.4) {\footnotesize $\calC_{10}$ $(P=0)$};
    \draw[->] (0,-4) -- (0,4) node[above left] {\footnotesize $P$};
    \node at(0,-4.4) {\footnotesize $\calC_{11}$ $(S=0)$};
    \draw[densely dashed] plot (-\x,-\x) node[below left] {\footnotesize $\calC_{\min}^\perp$ $(P=S)$};
    \draw[dotted] plot (\x,-\x) node[below right] {\footnotesize $P=-S$};
    \draw[dotted] plot[domain=-2:2] (\x,2*\x) node[above right] {\footnotesize anti-quasi-Sasakian $(P=2S)$};
    \draw[densely dashed] plot[domain=-2:2] (-\x,2*\x) node[above left] {\footnotesize $\calC_{\min}$ $(P=-2S)$};
\end{tikzpicture}
\caption{Remarkable classes of type $\calC_{\lambda,\mu}$}
\label{fig:C_lm}
\end{figure}

%----------------------------------
\subsubsection{Structures with integrable horizontal distribution} \label{deta=0}

We now consider $\calH$-parallel structure for which the horizontal distribution $\calH$ is integrable, or equivalently, the $1$-form $\eta$ is closed. This case is characterized by the following:

\begin{proposition}
	Let $(M,\f,\xi,\eta,g)$ be an $\calH$-parallel almost contact metric manifold. Then the following conditions are equivalent:		
	\begin{enumerate}
		\item $d\eta=0$;
		\item $S\xi=0$ and $\f S\in S^2\calH^*$ (or equivalently $\f S\in S^2TM^*$);
		\item $S\xi=0$ and $S\f \in S^2\calH^*$;
		\item $M$ is of class $\calC_5\oplus\calC_8\oplus\calC_{9}\oplus\calC_{11}$.
	\end{enumerate}
	In particular, $\calC_{11}$ coincides with the class of $\calH$-parallel almost contact metric structures, with $\xi$ Killing and $d\eta=0$.
    All in all, the class of almost contact metric manifolds with $d\eta=0$ is $\calD_1\oplus\calC_5\oplus\calC_8\oplus\calC_{9}\oplus\calC_{11}$.
\end{proposition}
\begin{proof}
	The equivalence between 1 and 2, immediately follows from the expression of $d\eta$ written in \eqref{eq:deta}. Applying it to $\f X$ and $\f Y$, and using $\eta\circ S=0$, we have
	$$d\eta(\f X,\f Y)=-g(S\f X,Y)+g(X,S\f Y),$$
	and then $d\eta$ vanishes on $\calH$ if and only if $S\f$ is symmetric with respect to $g$. Moreover, by Proposition \ref{Prop:char.D2},  $d\eta(\xi,\cdot)=0$ if and only if $S\xi=0$. As regards the equivalence of 3 and 4, consider the space 
	$$\calD:=\{\alpha\in\calD_2\ |\ \f S\in S^2\calH^*\}\subset\calD_2.$$
	From the characterizations of $\calC_5$, $\calC_8$, $\calC_{9}$, $\calC_{11}$, it is clear that  $\calC_5\oplus\calC_8\oplus\calC_{9}\oplus\calC_{11}\subset \calD$. On the other hand, by Proposition \ref{Prop:Killing}, the intersection of $\calD$ with the class $\calC_6\oplus\calC_7\oplus\calC_{10}\oplus\calC_{11}$ is exactly $\calC_{11}$ ($S=0$). Therefore $\calD=\calC_5\oplus\calC_8\oplus\calC_{9}\oplus\calC_{11}$.
   The last statement is an immediate consequence of the fact that $d\eta=0$ for structures of class $\calD_1$ (Remark \ref{remark-d1}).
\end{proof}

%----------------------------------
\section{The intrinsic torsion}
%----------------------------------
\subsection{Intrinsic torsion and minimal connection of almost contact metric manifolds}

In this section we will provide a description of the intrinsic torsion of an $\calH$-parallel almost contact metric manifold via the intrinsic endomorphisms $(S,h)$. We will use the fact that the intrinsic torsion is encoded by the minimal connection introduced by J. C. Gonz\'alez-D\'avila and F. Mart\'in Cabrera in \cite{G.MC.}.\\

First we recall the definition of the intrinsic torsion of an almost contact metric manifold $M^{2n+1}$.
An almost contact metric structure $(\f,\xi,\eta,g)$ corresponds to a reduction of the structure group of the frame bundle to $G:=U(n)\times 1$, where $U(n)\times 1\subset SO(2n+1)$. The Lie algebra $\so(2n+1)$ can be decomposed into the sum of the Lie algebra $\fraku(n)$ of $G$ and its orthogonal complement $\m$, i.e. 
$$\so(2n+1)=\fraku(n)\oplus\m.$$ 
Identifying $\so(2n+1)\cong\Lambda^2M$, we have that
$$\fraku(n)\cong\{\omega\in\Lambda^2M\ |\ \omega(\f X,\f Y)=\omega(X,Y)\ \forall X,Y\in TM \}$$
is the subspace of $\f$-invariant $2$-forms, and
$$\m\cong\{\omega\in\Lambda^2M\ |\ \omega(\f X,\f Y)=-\omega(X,Y)\ \forall X,Y\in\calH\}$$
is the subspace of the $2$-forms which are $\f$-anti-invariant on $\calH$.
Taking an orthonormal frame adapted to the structure, i.e. a $\f$-basis $\{e_1,\dots,e_{2n+1}=\xi\}$, consider the $1$-form $\Omega^g$ with values in the Lie algebra $\so(2n+1)$ defined by the connection forms of the Levi-Civita connection, namely
$$\Omega^g:=(\omega^g_{ij}), \qquad \omega^g_{ij}:=g(\nablag e_i, e_j), \quad i,j=1,\dots,2n+1.$$
The \textit{intrinsic torsion} of the almost contact metric structure is defined as $$\Gamma:=\pr_\m(\Omega^g),$$
where $\pr_\m:\so(2n+1)\to\m$ denotes  the orthogonal projection onto $\m$.
Now, let us consider a connection $\nabla$ on $M^{2n+1}$, and denote by $A$ the difference tensor field between $\nabla$ and the Levi-Civita connection, i.e. $A:=\nabla-\nabla^g.$ We denote by the same symbol the $(0,3)$-tensor field defined by $A(X,Y,Z):=g(A(X,Y),Z)$. It is known that $\nabla$ is a metric connection ($\nabla g=0$) if and only if
\begin{equation}\label{eq:A_metric conn.}
	A(X,Y,Z)+A(X,Z,Y)=0.
\end{equation}

As for the Levi-Civita connection, one can consider the $1$-form with values in $\so(2n+1)$ given by
$$\Omega:=(\omega_{ij}),\qquad \omega_{ij}:=g(\nabla e_i,e_j).$$
Then, for every $X\in TM$, $$\Omega(X)=\Omega^g(X)+A(X),$$
where $A(X):=A(X,\cdot,\cdot)$ is a $2$-form because of \eqref{eq:A_metric conn.}.
The metric connection $\nabla$ is said to be \emph{adapted} (to the almost contact structure) if it parallelizes all the structure tensor fields, that is
\begin{equation}\label{eq:adapted}
	\nabla g=0,\quad \nabla\f=0,\quad \nabla\xi=0,\quad \nabla\eta=0.
\end{equation}
This is equivalent to requiring that the projection of $\Omega$ onto $\m$ vanishes, $$\pr_\m(\Omega(X))=\Gamma(X)+\pr_\m(A(X))=0,$$ that is
$\Gamma(X)=-\pr_\m(A(X))$ for every $X\in TM$.  
\\

In \cite{G.MC.} the authors introduce a special metric connection adapted to any almost contact metric structure which is called the \emph{minimal} connection. We recall the result, providing an explicit proof.

\begin{theorem}\label{Thm:connection}
	Let $(M,\f,\xi,\eta,g)$ be an almost contact metric manifold. Then there exists a unique adapted metric connection $\nabla=\nablag+A$,  such that
	\begin{equation}\label{eq:intrinsic_torsion}
		A(X,\f Y,\f Z)+A(X,Y,Z)=0,\quad \forall X\in TM,\ \forall Y,Z\in\calH
	\end{equation}
	or equivalently, for every $X\in TM$, $A(X)\in\m$. Consequently, the intrinsic torsion of the structure is $\Gamma(X)=-\pr_\m(A(X))=-A(X)$.
	The tensor field $A$ is given by
	\begin{align}
		2A(X,Y)&=-\f(\nablag_X\f)Y+2((\nablag_X\eta)Y)\xi-\eta(Y)\nablag_X\xi \label{eq:A_(1)}\\
		&=(\nablag_X\f)\f Y+((\nablag_X\eta)Y)\xi- 2\eta(Y)\nablag_X\xi.\label{eq:A_(2)}
	\end{align}
\end{theorem}
\begin{proof}
	One easily checks that the stated formula for $A$ does indeed define a connection with the desired properties. Conversely, assume that a connection $\nabla$ as in the statement exists. From $\nabla\xi=0$, we conclude $A(X,\xi)=-\nablag_X\xi$ for every $X\in TM$. Moreover, since $\nabla g=0$, \eqref{eq:A_metric conn.} holds. We decompose any $X,Y,Z\in TM$, into $Y=Y_\calH+\eta(Y)\xi$ and $Z=Z_\calH+\eta(Z)\xi$, so that \eqref{eq:intrinsic_torsion} gives
	\begin{align*}
		&A(X,\f Y,\f Z)+A(X,Y,Z)\\
		&=A(X,\f Y_\calH,\f Z_\calH)+A(X,Y_\calH,Z_\calH)+ \eta(Y)A(X,\xi,Z)+\eta(Z)A(X,Y,\xi)\\
		&=\eta(Y)g(A(X,\xi),Z)-\eta(Z)g(A(X,\xi),Y)\\ 
		&={}-\eta(Y)g(\nablag_X\xi,Z)+\eta(Z)g(\nablag_X\xi,Y),
	\end{align*}
	which implies the identity
	\begin{equation}\label{eq:A-fA}
		A(X,Y)-\f A(X,\f Y)=-\eta(Y)\nablag_X\xi+g(\nablag_X\xi,Y)\xi.
	\end{equation}
	On the other hand, $\nabla\f=0$ if and only if
	$$(\nablag_X\f)Y=\f A(X,Y)-A(X,\f Y).$$
	Applying $\f$ to both sides of this equation, we get
	\begin{align*}
		\f(\nablag_X\f)Y&=-A(X,Y)+\eta(A(X,Y))\xi-\f A(X,\f Y)\\
		&=-A(X,Y)+A(X,Y,\xi)\xi-\f A(X,\f Y)\\
		&=-A(X,Y)+g(\nablag_X\xi,Y)\xi-\f A(X,\f Y),
	\end{align*}
	that is
	$$A(X,Y)+\f A(X,\f Y)=\f(\nablag_X\f)Y+g(\nablag_X\xi,Y)\xi.$$
	Comparing with \eqref{eq:A-fA}, and using $g(\nablag_X\xi,Y)=(\nablag_X\eta)Y$ one gets \eqref{eq:A_(1)}. Finally, \eqref{eq:A_(2)} follows from
	\[(\nablag_X\f)\f Y+\f(\nablag_X\f)Y= (\nablag_X\f^2)Y=((\nablag_X\eta)Y)\xi+\eta(Y)\nablag_X\xi.\qedhere\]
\end{proof}
\bigskip

Recall that the space $\mathcal{A}$ of all (0,3)-tensor fields satisfying \eqref{eq:A_metric conn.} splits into the orthogonal sum of three irreducible components under the action of the orthogonal group \cite{Agricola_Srni},
$$\mathcal{A}=\mathcal{A}_1\oplus\mathcal{A}_2\oplus\mathcal{A}_3,$$
where
\begin{align*}
	\mathcal{A}_1&=\{A\in\mathcal{A} \ |\ \exists\vartheta\in\Lambda^1(M)\ \text{s.t.}\ A(X,Y,Z)=g(X,Y)\vartheta(Z)-g(X,Z)\vartheta(Y)\}\cong TM,\\
	\mathcal{A}_2&=\{A\in\mathcal{A} \ |\ \mathfrak{S}_{X,Y,Z}A(X,Y,Z)=0,\; c_{12}(A)=0 \},\\
	\mathcal{A}_3&=\{A\in\mathcal{A} \ |\ A(X,Y,Z)+A(Y,X,Z)=0\} \cong\Lambda^3(M).
\end{align*}
Metric connections having $A\in \mathcal{A}_1$ are called connections of \textit{vectorial type}, those with $A\in\mathcal{A}_2$ are called connections of
\textit{Cartan type},
while metric connections with $A\in\mathcal{A}_3$ are called metric connection with \textit{totally skew-symmetric torsion} (or shortly \textit{skew torsion}).
\\
In this regard, recall that the spaces $\mathcal{A}$ is isomorphic to the space $\mathcal{T}$ of the torsion tensors $T$ of all possible metric connections, by means of
\begin{align}
	T(X,Y,Z)&=A(X,Y,Z)-A(Y,X,Z),\label{eq:T-A}\\
	2A(X,Y,Z)&=T(X,Y,Z)-T(Y,Z,X)+T(Z,X,Y),
\end{align}
where, as usual, $T(X,Y,Z):=g(T(X,Y),Z)$. 

\begin{remark}
	Notice that the $(0,3)$-tensor field $A$ satisfies the same symmetries of $\nablag\Phi$. In other words, the space $\Lambda^1(M)\otimes\m$ of possible intrinsic torsions coincides with $\calC(TM)$, for which we already have a decomposition into irreducible $U(n)\times\{1\}$-components. More precisely, we can further decompose $\m=\m_1\oplus\m_2$, with
	$$\m_1=\{\omega\in\Lambda^2M\ |\ \omega(\f X,\f Y)= -\omega(X,Y)\ \forall X,Y\in TM\},$$
	$$\m_2=\{\omega\in\Lambda^2M\ |\ \omega(X,Y)=0\ \forall X,Y\in\calH\}.$$ 
	It follows that
	\begin{align*}
		\calC(TM)&\cong\Lambda^1(M)\otimes\m=(\langle\eta\rangle\oplus\langle\eta\rangle^\perp)\otimes(\m_1\oplus\m_2)\\
		 &=(\langle\eta\rangle\otimes\m_1)\oplus(\langle\eta\rangle\otimes\m_2)\oplus(\langle\eta\rangle^\perp\otimes\m_1)\oplus(\langle\eta\rangle^\perp\otimes\m_2),
	\end{align*}
	where $\Lambda^1(M)=\langle\eta\rangle\oplus\langle\eta\rangle^\perp$ is the dual decomposition of $TM=\langle\xi\rangle\oplus\calH$, and $\langle\eta\rangle^\perp$ is the space of 1-forms vanishing on $\xi$. Then, taking into account Table \ref{Table:D-classes}, one realizes that:
	$$\calD_1\cong\langle\eta\rangle^\perp\otimes\m_1,\quad
	\calD_2\cong(\langle\eta\rangle^\perp\otimes\m_2)\oplus(\langle\eta\rangle\otimes\m_1),\quad
	\calD_3\cong\calC_{12}=\langle\eta\rangle\otimes\m_2.$$
	In particular, $\langle\eta\rangle^\perp\otimes\m_2=\calC_5\oplus\cdots\oplus\calC_{10}$ and $\langle\eta\rangle\otimes\m_1=\calC_{11}$.
\end{remark}

%-----------------------------------
\subsection{The torsion type of $\calH$-parallel structures}

We will now determine the intrinsic torsion type for each irreducible class of $\calH$-parallel almost contact metric manifolds. To this aim, first we specialize the formula \eqref{eq:A_(2)} in $\calH$-parallel case.

\begin{lemma}
Let $(M,\f,\xi,\eta,g)$ be an $\calH$-parallel almost contact metric manifold. Let $\nabla=\nabla^g+A$ be the minimal connection. Then, the tensor field $A$ is given by:
\begin{equation}\label{eq:intr.tors.Hp}
	2A(X,Y)=\eta(X)\f PY-\eta(X)\eta(Y)\f S\xi+2\eta(Y)\f SX+2g(SX,\f Y)\xi,
\end{equation}
or equivalently,
\begin{equation}\label{eq:intr.tors.Hp(0,3)}
	2A(X,Y,Z)=-\eta(X)g(PY,\f Z)+\eta(X)\eta(Y)g(S\xi,\f Z)+2\eta(Z)g(SX,\f Y) -2\eta(Y)g(SX,\f Z).
\end{equation}
\end{lemma}
\begin{proof}
	Taking the scalar product of \eqref{eq:A_(2)} with $Z\in TM$, using equations \eqref{eq:nablaPhi(1)}, \eqref{eq:S_nabla_eta}, $\nablag_X\xi=-\f SX$ and $P\xi=S\xi$, we get:
	\begin{align*}
		2A(X,Y,Z)&=2g(A(X,Y),Z)\nonumber\\
		&=-(\nablag_X\Phi)(\f Y,Z)+\eta(Z)(\nablag_X\eta)Y -2\eta(Y)g(\nablag_X\xi,Z)\nonumber\\
		&=\eta(X)g(P\f^2Y,\f Z)+\eta(Z)g(SX,\f Y)+\eta(Z)g(SX,\f Y)
		-2\eta(Y)g(SX,\f Z)\nonumber\\
		&=-\eta(X)g(PY,\f Z)+\eta(X)\eta(Y)g(S\xi,\f Z)+2\eta(Z)g(SX,\f Y)-2\eta(Y)g(SX,\f Z),
	\end{align*}
	which is also equivalent to \eqref{eq:intr.tors.Hp}.
\end{proof}

Taking into account the characterizing conditions of the subclasses of $\calD_2\oplus\calD_3$, we obtain the following:

\begin{theorem}\label{theorem-intrinsic}
	The irreducible classes of $\calH$-parallel almost contact metric manifolds correspond to $\nabla$ in:
	\begin{table}[H]
		\centering
		\begin{tabular}{c|c|c|c|c|c|c|c|c}
			\toprule
			Class of $M$& $\calC_5$&$\calC_6$&$\calC_7$&$\calC_8$&$\calC_9$&$\calC_{10}$&$\calC_{11}$&$\calC_{12}$\\\hline
			Type of $\nabla$&$\mathcal{A}_1$&$\mathcal{A}_2\oplus\mathcal{A}_3$&$\mathcal{A}_2\oplus\mathcal{A}_3$&$\mathcal{A}_2$&$\mathcal{A}_2$&$\mathcal{A}_2\oplus\mathcal{A}_3$&$\mathcal{A}_2\oplus\mathcal{A}_3$&$\mathcal{A}_1\oplus\mathcal{A}_2$\\
		\bottomrule
		\end{tabular}
	\end{table}
\end{theorem}
\begin{proof}
	In $\calC_5\oplus\cdots\oplus\calC_{10}$, being $S\xi=0$ and $P=0$, equation \eqref{eq:intr.tors.Hp(0,3)} reduces to
	\begin{equation}\label{eq:A_C5-C10}
		A(X,Y,Z)=\eta(Z)g(SX,\f Y)-\eta(Y)g(SX,\f Z).
	\end{equation}
	In particular, in $\calC_5$ we have that $S=\beta\f$, and hence
	\begin{align*}
		A(X,Y,Z)&=\beta[\eta(Z)g(\f X,\f Y)-\eta(Y)g(\f X,\f Z)]\\
		&=\beta[\eta(Z)g(X,Y)-\eta(Y)g(X,Z)],
	\end{align*}
	which says that $A\in \calA_1$, with $\vartheta=\beta\eta$.\\
	In $\calC_6$, $S=\alpha(I-\eta\otimes\xi)$. Therefore:
	\begin{align*}
		A(X,Y,Z)&=\alpha[\eta(Z)g(X,\f Y)-\eta(Y)g(X,\f Z)]\\
		&=\alpha[\eta(Z)\Phi(X,Y)+\eta(Y)\Phi(Z,X)]\\
		&=\alpha[(\eta\wedge\Phi)(X,Y,Z)-\eta(X)\Phi(Y,Z)],
	\end{align*}
	that is $A=\alpha(\eta\wedge\Phi-\eta\otimes\Phi)$.
	On the other hand, we can also write
	\begin{align*}
		A(X,Y,Z)&=\alpha[\eta(Z)\Phi(X,Y)+\eta(Y)\Phi(Z,X)-2\eta(X)\Phi(Y,Z)]+2\alpha\eta(X)\Phi(Y,Z)\\
		&=\alpha A'(X,Y,Z)+2\alpha\eta(X)\Phi(Y,Z),
	\end{align*}
	where we set $A'(X,Y,Z):=-2\eta(X)\Phi(Y,Z)+\eta(Y)\Phi(Z,X)+\eta(Z)\Phi(X,Y)$ for every $X,Y,Z\in TM$. Thus, $A=\alpha(A'+2\eta\otimes\Phi)$ and, comparing with the above expression, one gets that
	$$A=\frac13\alpha A'+\frac23\alpha\eta\wedge\Phi.$$
	Since $2\alpha\eta\wedge\Phi\in\Lambda^3(M)\cong\calA_3$, in order to prove that $A\in\calA_2\oplus\calA_3$, we only need to show that $A'\in\calA_2$. First we notice that it belongs to $\calA$, since it is skew-symmetric in the last two entries. Moreover:
	\begin{align*}
		\mathfrak{S}_{X,Y,Z}A'(X,Y,Z)&=
		-2\eta(X)\Phi(Y,Z)+\eta(Y)\Phi(Z,X)+\eta(Z)\Phi(X,Y)\\
		&\quad -2\eta(Y)\Phi(Z,X)+\eta(Z)\Phi(X,Y)+\eta(X)\Phi(Y,Z)\\
		&\quad -2\eta(Z)\Phi(X,Y)+\eta(X)\Phi(Y,Z)+\eta(Y)\Phi(Z,X)=0,
	\end{align*}
	and, with respect to a local orthonormal frame $\{e_i,\xi\}$ ($i=1,\dots,2n$),
	$$(c_{12}A')(Z)=\sum_{i=1}^{2n}A'(e_i,e_i,Z)+A'(\xi,\xi,Z)=\eta(Z)\sum_{i=1}^{2n}\Phi(e_i,e_i)=0.$$
	A similar argument holds for the classes $\calC_7$ and $\calC_{10}$. Indeed, in both cases $S\f$ turns out to be skew-symmetric with respect to $g$, so that $\Psi:=g(\cdot,S\f\cdot)$ defines a $2$-form such that $\Psi(\xi,\cdot)=0$ and
	$$A(X,Y,Z)=\pm[\eta(Z)\Psi(X,Y)+\eta(Y)\Psi(Z,X)],$$
	where the sign depends on $S\in S^2\calH^*$ or $S\in \Lambda^2\calH^*$.\\
	Now let us consider the case of $\calC_8$. Being $S\in\Lambda^2\calH^*$, \eqref{eq:A_C5-C10} can be written as
	$$A(X,Y,Z)=-\eta(Z)g(X,S\f Y)+\eta(Y)g(X,S\f Z).$$
	Being also $S\f=\f S$, it turns out that $S\f\in S^2\calH^*$ and hence $\mathfrak{S}_{X,Y,Z}A(X,Y,Z)=0$. Furthermore, fixed a local orthonormal frame $\{e_i,\xi\}$ ($i=1,\dots,2n$),
	$$(c_{12}A)(Z)=\sum_{i=1}^{2n}A(e_i,e_i,Z)=-\eta(Z)\sum_{i=1}^{2n}g(e_i,S\f e_i)=-\eta(Z)\tr(S\f)=0,$$
	being $\tr(\f S)=\tr(S\f)=0$ in $\calC_8$. Therefore, $A\in\calA_2$. \\
	The proof in the case of $\calC_9$ goes analogously. Indeed, even in this case $S\f$ is symmetric with respect to $g$, thus implying that the cyclic sum vanishes. Moreover, taking a local orthonormal frame of type $\{e_i,\f e_i,\xi\}$ ($i=1,\dots,n$), applying \eqref{eq:A_C5-C10} and $S\in S^2\calH^*$, we have:
	\begin{align*}
		(c_{12}A)(Z)&=\sum_{i=1}^{n}[A(e_i,e_i,Z)+A(\f e_i,\f e_i,Z)] =\eta(Z)\sum_{i=1}^{n}[g(Se_i,\f e_i)+g(S\f e_i,\f^2e_i)]\\
		&=\eta(Z)\sum_{i=1}^{n}[g(Se_i,\f e_i)-g(S\f e_i,e_i)]=0.
	\end{align*}
	Concerning $\calC_{11}$, which is characterized by $S=0$, we have that $P=-2h$ and then, by equation \eqref{eq:intr.tors.Hp(0,3)}, it turns out that
	$$A(X,Y,Z)=\eta(X)g(hY,\f Z)=\eta(X)g(h\f Y,Z),$$
	where we used the fact that $h$ anticommutes with $\f$. By the skew-symmetry of $A$ in $Y$ and $Z$, one deduces that $h\f\in\Lambda^2TM^*$, thus defining a $2$-form $\omega:=g(h\f\cdot,\cdot)$. Straightforward computations show that $3A=\eta\wedge\omega+A',$ where $A'(X,Y,Z)=2\eta(X)\omega(Y,Z)-\eta(Y)\omega(Z,X)-\eta(Z)\omega(X,Y)$ belongs to $\calA_2$, thus proving that $A\in\calA_2\oplus\calA_3$.\\
	Finally, in $\calC_{12}$, being $\f h=0$ and $S|_\calH=0$ (i.e. $S\f=0$), then $\f P=\f S=\eta\otimes\f S\xi$ and \eqref{eq:intr.tors.Hp(0,3)} reduces to
	$$A(X,Y,Z)=\eta(X)\eta(Y)g(\f S\xi,Z)-\eta(X)\eta(Z)g(\f S\xi,Y).$$
	Assuming $\dim M=2n+1$, let us define the $1$-form  $\vartheta=\frac1{2n}g(\f S\xi,\cdot)$, so that the $(0,3)$-tensor field $A_1(X,Y,Z):=g(X,Y)\vartheta(Z)-g(X,Z)\vartheta(Y)$ belongs to $\calA_1$. In order to prove the statement, it suffices to show that $A_2:=A-A_1$ belongs to $\calA_2$. For every $X,Y,Z\in TM$
	\begin{align*}
		A_2(X,Y,Z)&=A(X,Y,Z)-A_1(X,Y,Z)\\
		&=(2n\eta(X)\eta(Y)-g(X,Y))\vartheta(Z)-(2n\eta(X)\eta(Z)-g(X,Z)) \vartheta(Y).
	\end{align*}
	Then it is clear that $\mathfrak{S}_{X,Y,Z}A_2(X,Y,Z)=0$. Moreover,  fixed a local orthonormal frame $\{e_i,\xi\}$ ($i=1,\dots,2n$), since $\vartheta(\xi)=0$, one has:
	\begin{align*}
		(c_{12}A_2)(Z)&=\sum_{i=1}^{2n}A_2(e_i,e_i,Z)+A_2(\xi,\xi,Z) \\ &=-2n\vartheta(Z)+\sum_{i=1}^{2n}g(e_i,Z)\vartheta(e_i)
		+(2n-1)\vartheta(Z)=0. \qedhere
	\end{align*}
\end{proof}

\begin{remark}
	It is worth remarking that, except for $\calC_5$, $\calC_8$ and $\calC_9$, the intrinsic torsion corresponding to the other classes lies in a sum class, and it cannot be of pure type, unless the manifold is cok\"ahler (in which case $A=0$, i.e. $\nabla=\nablag$). This easily follows from the explicit decomposition of $A$ computed class by class in the above proof.
\end{remark}

%-------------------------------------
\subsection{The case of totally skew-symmetric torsion}

Our last goal in this paper is to investigate in more detail adapted connections with skew torsion, both from the point of view of minimal connections (as just described) and characteristic connections (to be described immediately); it will turn out that both are related in a very satisfying way and that
the subspace $\calC_{\min}\subset\calC_{10}\oplus\calC_{11}$ is the piece linking the two.

Recall that Friedrich and Ivanov provided necessary and sufficient conditions for any almost contact metric manifold $(M,\f,\xi,\eta,g)$ to admit a metric connection with skew-symmetric torsion preserving the almost contact structure (i.e. satisfying \eqref{eq:adapted}) in \cite{FI02}. Such a connection is then necessarily unique and called \textit{characteristic}, and precisely they proved the following:
\smallskip

\begin{theorem}[Friedrich, Ivanov]\label{Th:FI02}
	Let $(M,\f,\xi,\eta,g)$ be an almost contact metric manifold. It admits a characteristic connection $\nabla$ if and only if $\xi$ is a Killing vector field and $N_\f$ is totally skew-symmetric. The connection $\nabla$ is uniquely determined and its torsion is given by
	$$T=\eta\wedge d\eta+N_\f+d^\f\Phi-\eta\wedge(\xi\lrcorner N_\f),$$
	where $d^\f\Phi(X,Y,Z):=d\Phi(\f X,\f Y,\f Z)$ and $(\xi\lrcorner N_\f)(X,Y)=N_\f(\xi,X,Y)$.
\end{theorem}

Another equivalent condition for the existence of a characteristic connection on an almost contact metric manifold can be found in \cite{Agr.Holl}. Here the authors also discuss which Chinea-Gonzalez' classes may admit or not a characteristic connection.\\
Now we specialize these results in the case of $\calH$-parallel almost contact metric manifolds.

\begin{lemma}\label{Prop:char.connection}
	Let $(M,\f,\xi,\eta,g)$ be an $\calH$-parallel almost contact metric manifold. Then, $M$ admits a characteristic connection if and only if $\xi$ is Killing and
	\begin{equation}\label{eq:skew-symm.N}
		h\f= S\f-\f S, \ \text{or, equivalently, } \ P=-h.
	\end{equation} 	
\end{lemma}
\begin{proof} 
	If $\xi$ is Killing, the equivalence of the two equations in \eqref{eq:skew-symm.N} is an immediate consequence of the definition of $P$ and $S\xi=0$.	According to Theorem \ref{Th:FI02}, we need to show that, when $\xi$ is a Killing vector field, equation \eqref{eq:skew-symm.N} is equivalent to require that $N_\f$ is totally skew-symmetric. Since $N_\f$ is skew-symmetric as $(1,2)$-tensor field, $N_\f\in\Lambda^3(M)$ if and only if it is skew-symmetric in the last two entries. Furthermore, since $N_\f$ vanishes on triplets of horizontal vector fields (by Proposition \ref{Prop:dPhi=N=0}), we only need to show that \eqref{eq:skew-symm.N} is equivalent to $$N_\f(X,Y,\xi)=-N_\f(X,\xi,Y)\quad \forall X,Y\in TM.$$
	Being $\xi$ Killing, both $S\f$ and $\f S$ are skew-symmetric with respect to $g$ (see Proposition \ref{Prop:Killing}). Then, by equation \eqref{N1}, 
    $$N_\f(X,Y,\xi)=2g((S\f-\f S)X,Y),$$
	$\xi$ Killing also implies $d\eta(\xi,\cdot)=0$ and $\f h=-h\f$. Then, from \eqref{eq:N(xi,.)}, we have
	$$N_\f(X,\xi,Y)=-2g(h\f X,Y),$$
	Comparing the two obtained expressions, we get \eqref{eq:skew-symm.N}.
\end{proof}
\medskip

One could wonder whether the characteristic connection coincides with the minimal connection. Owing to the uniqueness of the characteristic connection, the following proposition clarifies this aspect. This is in fact where the subspace $\calC_{\min}$ defined before appears, and the result which explains its name:

\begin{proposition}\label{Prop.H-parallel_A3}
Let $(M,\f,\xi,\eta,g)$ be an $\calH$-parallel almost contact metric manifold. Then, the following conditions are equivalent:
	\begin{enumerate}
		\item the minimal connection is of type $\mathcal{A}_3$;
		\item $M$ belongs to the class $\calC_{\min} \subset \calC_{10}\oplus\calC_{11}$.
	\end{enumerate}
\end{proposition}
\begin{proof}
	Assume that the minimal connection is of type $\mathcal{A}_3$. Then, by Proposition \ref{Prop:char.connection}, $\xi$ is Killing and $h\varphi=S\varphi-\varphi S$, or equivalently $P=-h$. In particular, $M$ belongs to the class $\calC_6\oplus\calC_7\oplus\calC_{10}\oplus\calC_{11}$. Now, from \eqref{eq:intr.tors.Hp(0,3)}, for every $X,Z\in\calH$, one has
	$$2A(X,\xi,Z)=-2g(SX,\varphi X),\qquad 2A(\xi,X,Z)=-g(PX,\varphi Z),$$
	so that $P=-2S$, and thus $h=2S$. Notice that $S\in\Lambda^2\calH^*$ since $P$ is skew-symmetric, and $S\varphi+\varphi S=0$, since $h$ anti-commutes with $\varphi$. By Proposition \ref{Prop:C10+C11}, $M$ belongs to $\calC_{10}\oplus\calC_{11}$.

	Conversely, assume that $M$ belongs to $\calC_{10}\oplus\calC_{11}$ and $h=2S$. Then $S\in\Lambda^2\calH^*$ and $S\varphi+\varphi S=0$. By the definition of $P$, we have
	$$P=-2h+S+\varphi S\varphi=-2h+2S=-2S.$$
	A direct computation applying \eqref{eq:intr.tors.Hp(0,3)} shows that
	$$A(X,Y,Z)=-\eta(X)g(\varphi SY,Z)-\eta(Y)g(\varphi SZ,X)-\eta(Z)g(\varphi SX,Y),$$
	which is a $3$-form since $\varphi S$ is skew-symmetric. Notice that, being $d\eta(X,Y)=-2 g(\varphi SX,Y)$, the skew-torsion of the minimal connection is $T=2A=\eta\wedge d\eta$, which we will see to be coherent with the next theorem.
\end{proof}

\begin{remark}
	In \cite{Pu13}, C. Puhle studied a special class of $5$-dimensional almost contact metric manifolds denoted by $\mathcal{W}_4\oplus\mathcal{W}_7$, which coincides with the Chinea-Gonzalez class $\calC_{10}\oplus\calC_{11}$, as stated in \cite[Theorem 3.1]{Pu12}. Within $\mathcal{W}_4\oplus\mathcal{W}_7$, the subclass $\mathcal{W}_4$ is characterized by
	$$d\Phi(\xi,\f X,Y)=-2d\eta(X,Y)$$
	(see \cite[Theorem 4.1]{Pu13}). Now, taking into account that $\xi$ is Killing, by using \eqref{eq:dPhi-Lg-h} and \eqref{eq:deta}, one gets
	$$g(\f X,hY)=-g(SX,\f Y)+g(\f X,SY).$$
	Applying Proposition \ref{Prop:C10+C11}, this is equivalent to $h=2S$. Therefore, according to Proposition \ref{Prop.H-parallel_A3}, the class of $\calH$-parallel almost contact metric $5$-manifolds in $\calC_{\min}$  coincides with $\mathcal{W}_4$.
\end{remark}
\medskip

The class $\calC_{\min}$ allows us to formulate a final characterization in terms of Chinea-Gonzalez classes for the existence of the characteristic connection:

\begin{theorem}[Existence of the characteristic connection] \label{new-char-conn}
\begin{enumerate}
    \item[] 
    \item An $\calH$-parallel almost contact metric manifold admits a characteristic connection $\nabla$ if and only if it is of class $\calC_6\oplus\calC_7\oplus\calC_{\min}$, and furthermore the torsion is then given by the expression $T=\eta\wedge d\eta$.
    \item An almost contact metric manifold admits a characteristic connection $\nabla$ if and only if it is of class $\calC_1\oplus\calC_3\oplus\calC_4\oplus\calC_6\oplus\calC_7\oplus\calC_{\min}$.
\end{enumerate}
\end{theorem}
\begin{proof}
	Let us denote by $\mathcal{W}$ the class of $\calH$-parallel almost contact metric manifolds admitting the characteristic connection. Using Lemma \ref{Prop:char.connection}, the characterization of the $\xi$ Killing condition, and that fact that both $P$ and $h$ vanish on $\calC_6\oplus\calC_7$, we have
	$$\calC_6\oplus\calC_7\subset \mathcal{W}\subset\calC_6\oplus\calC_7\oplus\calC_{10}\oplus\calC_{11}.$$
	Now, take $\alpha=\alpha_1+\alpha_2$, with $\alpha_1\in\calC_6\oplus\calC_7$ and $\alpha_2\in\calC_{10}\oplus\calC_{11}$. We show that $\alpha\in\mathcal{W}$ if and only if $\alpha_2\in\calC_{\min}$. Considering the associated operators $S$, $P$, $h$, and $S_i$, $P_i$, $h_i$, $i=1,2$, in this case we have $P_1=0$ and $h_1=0$, so that $P=P_2$ and $h=h_2$ (see Remark \ref{remark-sum}). By Lemma \ref{Prop:char.connection}, $\alpha\in\mathcal{W}$ if and only if $P=-h$, which means $P_2=-h_2$. In the class $\calC_{10}\oplus\calC_{11}$ this is equivalent to $P_2=-2S_2$, namely $\alpha_2\in\calC_{\min}$.

	Assuming that $M$ admits a characteristic connection, we determine its torsion. From \eqref{eq:dPhi=N=0} it follows that $d^\f\Phi=0$ and
	\begin{align*}
		N_\f(X,Y,Z)&=\eta(X)N_\f(\xi,Y,Z)+\eta(Y)N_\f(X,\xi,Z)+\eta(Z)N_\f(X,Y,\xi)\\
		&=\eta(X)N_\f(\xi,Y,Z)+\eta(Y)N_\f(\xi,Z,X)+\eta(Z)N_\f(\xi,X,Y)\\
		&=(\eta\wedge(\xi\lrcorner N_\f))(X,Y,Z),
	\end{align*}
	where we used the fact that $N_\f\in\Lambda^3(M)$. Hence, by Theorem \ref{Th:FI02}, we have that $T=\eta\wedge d\eta$.

    Let us now consider the general case of an arbitrary almost contact metric manifold. Assume $\xi$ Killing and let $\alpha=\alpha_1+\alpha_2$ be the $(0,3)$-tensor defined by the covariant derivative $\nabla^g\Phi$, where $\alpha_1\in\calD_1$ and $\alpha_2\in\calC_6\oplus\calC_7\oplus\calC_{10}\oplus\calC_{11}$. Now, $M$ admits a characteristic connection if and only if $N_\f$ is totally skew-symmetric.  On the one hand, $N_\f$ is totally skew-symmetric on $\calH$ if and only if the $\calD_1$-component $\alpha_1$ belongs to $\calC_1\oplus\calC_3\oplus\calC_4$ (see Remarks \ref{remark-total-class} and \ref{remark-d1}). On the other hand, by the previous discussion, the skew-symmetry $N_\f(X,Y,\xi)=-N_\f(X,\xi,Y)$ is equivalent to require that $\alpha_2$ belongs to $\calC_6\oplus\calC_7\oplus\calC_{\min}$ (see again Remark \ref{remark-total-class}). Observe that in this case, the formula for the torsion cannot be simplified any further, i.e.~it is given by the expression in Theorem \ref{Th:FI02}.
\end{proof}

We now characterize the $\nabla$-parallelism of the torsion $T=\eta\wedge d\eta$ for any $\calH$-parallel almost contact metric manifold admitting a characteristic connection. This is given in terms of the Riemannian curvature $R^g$ of the manifold $(M,g)$.

\begin{theorem}\label{thm:parallel-torsion}
    Let $(M,\f,\xi,\eta,g)$ be an $\calH$-parallel almost contact metric manifold admitting a characteristic connection $\nabla$ (i.e. $M$ of class $\calC_6\oplus\calC_7\oplus\calC_{\min}$). Then, the following statements are equivalent:
    \begin{enumerate}[label=(\roman*)]
        \item $\nabla T=0$, 
        \item $R^g(\xi,X)Y = -\eta(Y)\tilde X + g(\tilde X,Y)\xi$,
        \item $R^g(X,Y)\xi=\eta(Y)\tilde X-\eta(X)\tilde Y$,
    \end{enumerate}
    for every $X,Y\in TM$, where $\tilde X:= (S^2-hS)X$. In particular, $S^2-hS=S^2$ or $S^2-hS=-S^2$ if $M$ belongs to $\calC_6\oplus\calC_7$ or $\calC_{\min}$, respectively.
\end{theorem}
\begin{proof}
    To make use of simpler notation, let us set $\psi:=-\f S=\nablag\xi$ for the rest of the proof. Taking into account \eqref{eq:skew-symm.N}, observe that
    \begin{equation}\label{eq:psi2}
        \psi^2=\f S\f S=(S\f-h\f)\f S=S\f^2 S-h\f^2S=-S^2+hS.
    \end{equation}
    In particular, in $\calC_6\oplus\calC_7$, $h=0$ gives $\psi^2=-S^2$, while in $\calC_{\min}$, $h=2S$ gives $\psi^2=S^2$, which justifies the last statement.
    
    Now, since $\xi$ is Killing, $\psi$ is skew-symmetric with respect to $g$, and $\psi^2$ is symmetric. Moreover, by \eqref{eq:deta}, 
    \begin{equation}\label{eq:deta_psi}
        d\eta(X,Y)=g(X,\f SY)-g(\f SX, Y)=-2g(\f SX,Y)=2g(\psi X,Y).
    \end{equation} 
    Being $\nabla\eta=0$ and $\nabla g=0$, the parallelism of $T=\eta\wedge d\eta$ with respect to $\nabla$ is equivalent to the parallelism of $d\eta$, which in turn is equivalent to the parallelism of $\psi$. Thus, we show that $\nabla\psi=0$ if and only if (ii) or (iii) holds.
    Now, owing to \eqref{eq:deta_psi}, the difference between $\nabla$ and $\nablag$ is given by
    \begin{align*}
        A(X,Y,Z)&=\frac12 T(X,Y,Z)=\frac12(\eta\wedge d\eta)(X,Y,Z)=\\
        &=\frac12[\eta(X)d\eta(Y,Z)+\eta(Y)d\eta(Z,X)+\eta(Z)d\eta(X,Y)]\\
        &=\eta(X)g(\psi Y,Z)+\eta(Y)g(\psi Z,X)+\eta(Z)g(\psi X,Y),
    \end{align*}
    which gives that
    $$A(X,Y)=\eta(X)\psi Y-\eta(Y)\psi X+g(\psi X,Y)\xi.$$
    Moreover, by a well-known property of Killing vector fields (see \cite[Chap. VI, Proposition 2.6]{KN1}), the following holds:
    \begin{equation}\label{eq:R(xi,X)Y}
        (\nablag_X\psi)Y=-R^g(\xi,X)Y.
    \end{equation} 
    Then, using $\nabla=\nablag+A$, for every $X,Y\in TM$ we have:
    \begin{align*}
        (\nabla_X\psi)Y&=(\nablag_X\psi)Y+A(X,\psi Y)-\psi A(X,Y)\\
        &=-R^g(\xi,X)Y+\eta(X)\psi^2Y+g(\psi X,\psi Y)\xi-\eta(X)\psi^2 Y+\eta(Y)\psi^2X\\
        &=-R^g(\xi,X)Y+\eta(Y)\psi^2X-g(\psi^2X,Y)\xi.
    \end{align*}
    Therefore, owing to \eqref{eq:psi2}, the vanishing of $\nabla\psi$ is equivalent to (ii). Applying the Bianchi identity, one also has that (ii) implies (iii). Indeed, 
    \begin{align*}
        R^g(X,Y)\xi&=-R^g(\xi,X)Y+R^g(\xi,Y)X\\
             &=\eta(Y)(S^2-hS)X-\eta(X)(S^2-hS)Y.
    \end{align*}
    Conversely, assuming (iii), taking the scalar product with $Z\in TM$, and applying the symmetries of the Riemannian curvature, one has:
    \begin{align*}
        &g(R^g(X,Y)\xi,Z)=\eta(Y)g((S^2-hS)X,Z)-\eta(X)g((S^2-hS)Y,Z)\\
        &\Leftrightarrow g(R^g(\xi,Z)X,Y)=\eta(Y)g((S^2-hS)X,Z)-\eta(X)g(Y,(S^2-hS)Z)
    \end{align*}
    for every $X,Y,Z\in TM$, namely
    $$R^g(\xi,Z)X=-\eta(X)(S^2-hS)Z+g((S^2-hS)X,Z)\xi,$$
    which is equal to (ii) with $Z,X$ in place of $X,Y$.
\end{proof}

\begin{remark}
	In the class $\calC_6$ of $\alpha$-Sasakian manifolds, $d\eta=2\alpha\Phi$ and $T=2\alpha \eta\wedge\Phi$. Since $\nabla$ is adapted to the $(\f,\xi,\eta,g)$-structures, the torsion $T$ is $\nabla$-parallel if and only if $\alpha$ is constant. By \cite[Lemma 3.1]{Marrero}, this is always true if $\dim M\ge 5$. On the other hand, in the class $\calC_6$,  $S=\alpha(I-\eta\otimes\xi)$ (see Table \ref{Table:good_classes}). 
    Then, $$S^2-hS=S^2=\alpha^2(I-\eta\otimes\xi),$$ so that condition (iii) in the above theorem becomes
    $$R^g(X,Y)\xi=\alpha^2(\eta(Y)X-\eta(X)Y).$$
    For $\alpha=1$, this is a well-known identity for Sasakian manifolds \cite[Proposition 7.3]{blair10}.    
\end{remark}
\medskip

Next we recall an example of quasi-Sasakian manifold whose characteristic connection has parallel skew torsion. 

\begin{example}
	Let $\mathfrak{h}_\lambda$ be a $(2n+1)$-dimensional Lie algebra spanned by $\{\tau_i,\tau_{n+i},\xi\}$, $i=1,\dots,n$, whose only eventually non zero commutators are given by
    $$[\tau_i,\tau_{n+i}]=2\lambda_i\xi,$$
    for some $\lambda=(\lambda_1,\dots,\lambda_n)\in\R^n$. The corresponding simply connected Lie group $H^{2n+1}_\lambda$ is endowed with a left-invariant quasi-Sasakian structure $(\f,\xi,\eta,g)$, where $\f$ is defined by $$\f\tau_i=\tau_{n+i},\quad \f\tau_{n+i}=-\tau_i,\quad \f\xi=0,$$
    $g$ is the left-invariant metric with respect to which $\{\tau_i,\tau_{n+i},\xi\}$ is an orthonormal frame, and $\eta$ is the dual form of $\xi$. As stated in \cite[Section 9.3]{AFF}, the characteristic connection has parallel torsion.
   
    We indeed verify that condition (iii) of Theorem \ref{thm:parallel-torsion} holds.  Denoting by $\theta_i,\theta_{n+i}$ the dual forms of $\tau_i,\tau_{n+i}$, one has   $$d\eta=-2\sum_{i=1}^n\lambda_i\theta_i\wedge\theta_{n+i}.$$
    Taking into account \eqref{eq:deta_psi}, 
    $$\psi=\nablag\xi=-\sum_{i=1}^n\lambda_i(\theta_i\otimes\tau_{n+i}-\theta_{n+i}\otimes\tau_i),\qquad S^2=-\psi^2=\sum_{i=1}^n\lambda_i^2(\theta_i\otimes\tau_i+\theta_{n+i}\otimes\tau_{n+i}).$$
    Then,
    $$R^g(X,Y)\xi=\nablag_X\psi Y-\nablag_Y\psi X-\psi[X,Y]=\nablag_X\psi Y-\nablag_Y\psi X,$$
    where we used that the bracket is proportional to $\xi$ and $\psi\xi=\nablag_\xi\xi=0$.
   	A direct computation gives
   	$$\begin{array}{lll}
   	\nablag_{\tau_i}\tau_j=0, & \nablag_{\tau_i}\tau_{n+j}=\delta_{ij}\lambda_i\xi,&
  	\nablag_{\tau_i}\xi=\psi\tau_i=-\lambda_i\tau_{n+i},\\
   	\nablag_{\tau_{n+i}}\tau_j=-\delta_{ij}\lambda_i\xi, & \nablag_{\tau_{n+i}}\tau_{n+j}=0, &
   	\nablag_{\tau_{n+i}}\xi=\psi\tau_{n+i}=\lambda_i\tau_i,\\
   	\nablag_\xi\tau_i=-\lambda_i\tau_{n+i}&
   	\nablag_\xi\tau_{n+i}=\lambda_i\tau_i&
   	\nablag_\xi\xi=\psi\xi=0.
   	\end{array}$$
   	Therefore, for every $X,Y$ in the orthonormal frame,
   	$$R^g(X,Y)\xi=\eta(Y)S^2X-\eta(X)S^2Y.$$
\end{example}

%-------------------------------------

\bigskip\bigskip

%----------------------------------------------
{\sc Ilka Agricola}\\
Fachbereich Mathematik und Informatik, Philipps-Universit\"at Marburg, Hans-Meerwein Stra\ss e 6 (Campus Lahnberge), 35032 Marburg, Germany\\
\texttt{agricola@mathematik.uni-marburg.de}\\

{\sc Dario Di Pinto}\\
Dipartimento di Matematica, Università degli Studi di Bari Aldo Moro, Via E. Orabona 4, 70125 Bari, Italy\\
\texttt{dario.dipinto@uniba.it}\\

{\sc Giulia Dileo}\\
Dipartimento di Matematica, Università degli Studi di Bari Aldo Moro, Via E. Orabona 4, 70125 Bari, Italy\\
\texttt{giulia.dileo@uniba.it}\\

{\sc Marius Kuhrt} \\
Fachbereich Mathematik und Informatik, Philipps-Universit\"at Marburg, Hans-Meerwein Stra\ss e 6 (Campus Lahnberge), 35032 Marburg, Germany\\
\texttt{kuhrt@mathematik.uni-marburg.de}

\end{document}